\documentclass[11pt, twoside, final]{scrartcl}
\usepackage{amssymb}
\usepackage{amsfonts}
\usepackage{amsmath}
\usepackage{theorem}
\usepackage{pgfplots}
\usepackage{hyperref} \hypersetup{plainpages=false, colorlinks, linkcolor=black, citecolor=black, urlcolor=blue,
pdftitle={NFFT error},
pdfauthor={Daniel Potts, Manfred Tasche},
pdfstartview={FitBH}}

\usepackage{todonotes}

\theoremstyle{plain}
\newtheorem{theorem}{Theorem}[section]
\newtheorem{corollary}[theorem]{Corollary}
\newtheorem{lemma}[theorem]{Lemma}
\newtheorem{algorithm}[theorem]{Algorithm}
\newtheorem{definition}[theorem]{Definition}
\newtheorem{example}[theorem]{Example}
\newtheorem{remark}[theorem]{Remark}

\newcommand{\bend}{\hspace*{0ex} \hfill \hbox{\vrule height
    1.5ex\vbox{\hrule width 1.4ex \vskip 1.4ex\hrule  width 1.4ex}\vrule
    height 1.5ex}}
\newcommand{\qedsymbol}{\rule{1.5ex}{1.5ex}}

\newenvironment{Lemma}{\goodbreak\begin{lemma}\rmfamily}{\end{lemma}}
\newenvironment{Theorem}{\goodbreak\begin{theorem}}{\end{theorem}}
\newenvironment{Remark}{\goodbreak\begin{remark}\normalfont}{\bend\end{remark}}

\newenvironment{Example}{\goodbreak\begin{example}\normalfont}{\bend\end{example}}

\numberwithin{equation}{section}
\numberwithin{table}{section}
\numberwithin{figure}{section}

\renewcommand{\mathbf}[1]{\ensuremath{\boldsymbol{#1}}}

\title{Continuous window functions for NFFT}
\author{
  Daniel Potts\footnotemark[1]\and
  Manfred Tasche\footnotemark[3]
}

\date{}
\begin{document}
\maketitle

\begin{center}
\emph{Dedicated to J\"{u}rgen Prestin on the occasion of his 60th birthday}
\end{center}

\begin{abstract}
In this paper, we study the error behavior of the nonequispaced fast Fourier transform (NFFT). This approximate algorithm is mainly based on the convenient choice of a compactly supported window function. Here we consider
the continuous Kaiser--Bessel, continuous $\exp$-type, $\sinh$-type, and continuous $\cosh$-type window functions with the same support and same shape parameter.
We present novel explicit error estimates for NFFT with such a window function and derive rules for the optimal choice of the parameters involved in NFFT. The error constant of a window function depends mainly on the oversampling factor and the truncation parameter. For the considered continuous window functions, the error constants have an exponential
decay with respect to the truncation parameter.
\medskip

\smallskip

AMS \emph{Subject Classifications}:
65T50,  94A12,  42A10.   \end{abstract}

\footnotetext[1]{potts@mathematik.tu-chemnitz.de, Chemnitz University of
    Technology, Department of Mathematics, D--09107 Chemnitz, Germany}
\footnotetext[3]{manfred.tasche@uni-rostock.de, University of Rostock, Institute of Mathematics, D--18051 Rostock, Germany}

\section{Introduction}
The \emph{nonequispaced fast Fourier transform} (NFFT), see \cite{DuRo93, Bey95, St97, postta01} and \cite[Chapter 7]{PlPoStTa18} is an important generalization of the \emph{fast Fourier transform} (FFT). The window-based approximation leads to the most efficient algorithms under different approaches \cite{duro95,RuTo18}.
Recently a new class of window functions were suggested in \cite{BaMaKl18} and asymptotic error estimates are given in \cite{Ba20}. After \cite{St97}, the similarities of the window-based algorithms for NFFT became clear. Recently we have analyzed the window-based NFFT used so far and presented the related error estimates in \cite{PoTa20}.
Now we continue this investigation and present new error estimates for the some other window functions. More precisely, we consider the continuous Kaiser--Bessel window function
and two close relatives of the $\sinh$-type window function, namely the continuous $\exp$-type and $\cosh$-type window functions.
All these window functions have the same support and the same shape parameter.
We show that these window functions are very useful for NFFT, since they produce
very small errors.

In this paper, we present novel explicit error estimates \eqref{eq:estimateNFFT} with so-called error constants \eqref{eq:convenient1}. The error constants of NFFT are defined by values of the Fourier transform of the window function. We show that an upper
bound of \eqref{eq:convenient1} depends only on the oversampling factor $\sigma > 1$
and the truncation parameter $m \ge 2$ and decreases with exponential rate with respect to $m$. In numerous applications of NFFT, one uses quite often an oversampling factor $\sigma \in \big[\frac{5}{4},\,2\big]$ and a
truncation parameter $m \in \{2,\,3,\,\ldots,\,6\}$. Therefore we will assume that $\sigma \ge \frac{5}{4}$.

The outline of the paper is as follows. In Section \ref{Sec:Window} we introduce the set $\Phi_{m,N_1}$ of continuous, even window functions with support $\big[-\frac{m}{N_1},\,\frac{m}{N_1}\big]$,
where $m \in \mathbb N \setminus \{1\}$ and $N_1 = \sigma N\in 2\mathbb N$ (with $N \in 2 \mathbb N$ and $2m \ll N_1$) are fixed. We emphasize that a continuous window function $\varphi \in \Phi_{m,N_1}$
tends to zero at the endpoints $\pm \frac{m}{N_1}$ of its compact support $\big[-\frac{m}{N_1},\,\frac{m}{N_1}\big] \subset \big[-\frac{1}{2},\,\frac{1}{2}\big]$.
In Section \ref{Sec:rectWindow} we show that the simple rectangular window function \eqref{eq:rectwindow} is not convenient for NFFT.

The main results of this paper are contained in Sections \ref{Sec:KBWindow} -- \ref{Sec:coshWindow}.
For the first time, we present explicit estimates of the error constants \eqref{eq:convenient1} for fixed truncation parameter $m$
and oversampling factor $\sigma$.
In Section \ref{Sec:KBWindow}, we derive explicit error estimates for the continuous Kaiser-Bessel window function \eqref{eq:Kaiser-Bessel}. In comparison, we show that the popular standard
Kaiser--Bessel window function \eqref{eq:dKaiser-Bessel} has a similar error behavior as \eqref{eq:Kaiser-Bessel}.
A very useful continuous window function is the $\sinh$-type window function \eqref{eq:sinhwindow} which is handled in Section \ref{Sec:sinhWindow}.

The main drawback for the numerical analysis of the $\exp$-type and $\cosh$-type window function is the fact that an explicit
Fourier transform of this window function is unknown. In Sections \ref{Sec:expWindow} and \ref{Sec:coshWindow}, we develop a new technique. We split the continuous $\exp$-type/$\cosh$-type window function
into a sum $\psi + \rho$, where the Fourier transform of the compactly supported function $\psi$ is explicitly known and where the compactly supported function $\rho$ has small magnitude. Here we use the fact
that both window functions \eqref{eq:expwindow} and \eqref{eq:coshwindow} are close relatives of the $\sinh$-type window function \eqref{eq:sinhwindow} which was introduced by the authors in \cite{PoTa20}. The Fourier transform of $\rho$ is explicitly estimated for small as well as large frequencies, where $\sigma$ and $m$ are fixed. We present many numerical results so that the error constants of the different window functions can be easily compared. After this investigation, we favor the use of a continuous window function with small error constant, which can be very fast computed, such as the $\sinh$-type, standard/continuous $\exp$-type, or continuous $\cosh$-type window function.

\section{Continuous window functions for NFFT}\label{Sec:Window}

Let $\sigma > 1$ be an \emph{oversampling factor}.
Assume that $N\in 2 \,\mathbb N$ and $N_1:= \sigma N \in 2 \,\mathbb N$ are given. For fixed \emph{truncation parameter} $m \in {\mathbb N}\setminus \{1\}$
with $2m \ll N_1$, we introduce the open interval $I:=\big( - \frac{m}{N_1},\, \frac{m}{N_1}\big)$ and the set
$\Phi_{m,N_1}$ of all continuous window functions $\varphi:\, \mathbb R \to [0,\,1]$ with following properties:
\begin{quote}
$\bullet$ Each window function $\varphi$ is even, has the support $\bar I$,  and is continuous on~$\mathbb R$.\\
$\bullet$ Each restricted window function $\varphi|_{[0,\, m/N_1]}$ is decreasing with $\varphi(0) = 1$.\\
$\bullet$ For each window function $\varphi$, the Fourier transform
$$
{\hat \varphi}(v) := \int_{I} \varphi(t)\, {\mathrm e}^{-2\pi {\mathrm i}\,v t}\, {\mathrm d}t = 2\, \int_0^{m/N_1}\varphi(t)\, \cos(2\pi\,v t)\, {\mathrm d}t
$$
is positive for all $v\in [-N/2,\,N/2]$.
\end{quote}
Examples of continuous window functions of $\Phi_{m,N_1}$ are the (modified) B-spline window functions, algebraic window functions, Bessel window functions, and $\sinh$-type window functions (see \cite{PoTa20} and \cite[Chapter 7]{PlPoStTa18}). More examples are presented in Sections \ref{Sec:KBWindow} -- \ref{Sec:coshWindow}.
\medskip

In the following, we denote the torus $\mathbb R/\mathbb Z$ by $\mathbb T$ and the Banach space of continuous, 1-periodic functions by $C(\mathbb T)$. Let $I_N := \{-N/2, \ldots, N/2 -1 \}$ be the index set for
$N\in 2 \mathbb N$.

We say that a continuous window function $\varphi \in \Phi_{m,N_1}$ is \emph{convenient for} NFFT, if the $C(\mathbb T)$-{\emph{error constant}}
\begin{equation}
\label{eq:convenient1}
e_{\sigma}(\varphi) :=\sup_{N\in 2\mathbb N} e_{\sigma,N}(\varphi)
\end{equation}
with
\begin{equation}
\label{eq:esigmaNvarphi}
e_{\sigma,N}(\varphi) := \max_{n\in I_N}  \big\| \sum_{r\in \mathbb Z\setminus \{0\}} \frac{ {\hat \varphi}(n+r N_1)}{{\hat \varphi}(n)}\,{\mathrm e}^{2\pi {\mathrm i}\, r N_1\,\cdot} \big\|_{C(\mathbb T)}
\end{equation}
fulfills the condition $e_{\sigma}(\varphi) \ll 1$
for conveniently chosen oversampling factor $\sigma > 1$.

Now we show that the error of the nonequispaced fast Fourier transform (NFFT) with a window function $\varphi \in \Phi_{m,N_1}$ can be estimated by the error constant \eqref{eq:convenient1}. The NFFT (with nonequispaced spatial
data and equispaced frequencies) is an approximate, fast algorithm which computes approximately the values $p(x_j)$, $j=1,\ldots,M$, of any 1-periodic trigonometric
polynomial
$$
p(x) := \sum_{k \in I_N} c_k\, {\mathrm e}^{2 \pi {\mathrm i}\,k x}
$$
at finitely many nonequispaced nodes $x_j\in \big[-\frac{1}{2},\, \frac{1}{2}\big)$, $j=1,\ldots, M$, where $c_k \in \mathbb C$, $k\in I_N$, are given coefficients. By the properties of the window function $\varphi \in \Phi_{m,N_1}$, the 1-periodic function
$$
{\tilde \varphi}(x) := \sum_{k \in \mathbb Z} \varphi(x + k)\,, \quad x \in \mathbb T\,,
$$
is continuous on $\mathbb T$ and of bounded variation over $\big[-\frac{1}{2},\,\frac{1}{2}\big]$. Then from the Convergence Theorem of Dirichlet--Jordan (see \cite[Vol.~1, pp.~57--58]{Zy}), it follows that
$\tilde \varphi$ possesses the uniformly convergent Fourier expansion
$$
{\tilde \varphi}(x) = \sum_{k \in \mathbb Z} {\hat \varphi}(k)\,{\mathrm e}^{2 \pi {\mathrm i}\,k x}
$$
with the Fourier coefficients
$$
{\hat \varphi}(k) = \int_{\mathbb R} \varphi(x)\,{\mathrm e}^{-2 \pi {\mathrm i}\,k x}\,{\mathrm d}x = \int_0^1 {\tilde \varphi}(x)\,{\mathrm e}^{-2 \pi {\mathrm i}\,k x}\,{\mathrm d}x\,.
$$
We approximate
the trigonometric polynomial $p$  by the 1-periodic function
$$
s(x) := \sum_{\ell \in I_{N_1}} g_{\ell}\,{\tilde \varphi}\big(x - \frac{\ell}{N_1}\big) \in C(\mathbb T)
$$
with the coefficients
$$
g_{\ell} := \frac{1}{N_1} \sum_{k\in I_N} \frac{c_k}{{\hat \varphi}(k)}\, {\mathrm e}^{2 \pi {\mathrm i}\,k \ell/N_1}\,, \quad \ell \in I_{N_1}\,.
$$
The computation of the values $s(x_j)$, $j=1,\ldots, M$, is very easy, since $\varphi$ is compactly supported.
The computational cost of the algorithm is ${\cal O} (N_1\, \log N_1 + (2m+1)M)$, see \cite[Algorithm 7.1]{PlPoStTa18} and \cite{KeKuPo09} for details.

We interpret $s - p$ as the \emph{error function of the} NFFT which we measure in the norm of $C(\mathbb T)$.

\begin{Theorem}
\label{Thm:errorC+L2}
Let $\sigma > 1$, $N\in 2 \mathbb N$, and $N_1 = \sigma N \in 2 \mathbb N$ be given. Further let $m \in \mathbb N \setminus \{1\}$ with $2m \ll N_1$.

Then the error function of the $\mathrm{NFFT}$ can be estimated by
\begin{equation}
\| s - p\|_{C(\mathbb T)} \le e_{\sigma}(\varphi)\,\sum_{n \in I_N} |c_n|\,. \label{eq:estimateNFFT}
\end{equation}
\end{Theorem}

The proof of Theorem \ref{Thm:errorC+L2} is based on the equality
$$
s(x) - p(x) = \sum_{n \in I_N}\, \sum_{r \in \mathbb Z \setminus \{0\}} c_n\,\frac{{\hat \varphi}(n + r N_1)}{{\hat \varphi}(n)}\,{\mathrm e}^{2 \pi {\mathrm i}\,(n+r N_1)\,x}\,, \quad x \in \mathbb T\,.
$$
For details of the proof see \cite[Lemma 2.3]{PoTa20} and \cite[Lemma 2.1]{Ne14}.

In order to describe the behavior of $C(\mathbb T)$-error constant \eqref{eq:convenient1}, we have to study the Fourier transform of a window function $\varphi$ with the support $\bar I$.

\section{Rectangular window function}\label{Sec:rectWindow}

In this Section, we present a simple discontinuous window function which is not convenient for NFFT. Later we will use this discontinuous window function in Remarks \ref{Remark:KB-window} and \ref{Remark:origexpwindow}, where we estimate the $C(\mathbb T)$-error constants for the standard Kaiser-Bessel window function and the original $\exp$-type window function, respectively.

The simplest window function is the \emph{rectangular window function}
\begin{equation}
\label{eq:rectwindow}
\varphi_{\mathrm{rect}}(x) := \left\{ \begin{array}{ll} 1 &\quad x \in I\,,\\ [1ex]
\frac{1}{2} &\quad x = \pm \frac{m}{N_1}\,,\\ [1ex]
0 &\quad x \in \mathbb R \setminus \bar I\,,
\end{array}\right.
\end{equation}
Now we show that the rectangular window function \eqref{eq:rectwindow} is not convenient for NFFT.
The Fourier transform of \eqref{eq:rectwindow} has the form
\begin{equation}
\label{eq:FTrectwindow}
{\hat\varphi}_{\mathrm{rect}}(v) = \frac{2m}{N_1}\, \mathrm{sinc}\frac{2\pi m v}{N_1}\,, \quad v \in \mathbb R\,.
\end{equation}
The discontinuous window function ${\varphi}_{\mathrm{rect}}$ doesn't belong to $\Phi_{m,N_1}$.

\begin{Lemma}
\label{Lemma:FTrectwindow}
For each $n \in I_N \setminus \{0\}$, the Fourier series
$$
\frac{n}{N_1} \sum_{r\in \mathbb Z} \frac{1}{r + \frac{n}{N_1}} {\mathrm e}^{2\pi {\mathrm i}\, r N_1\,x}
$$
converges pointwise to the $\frac{1}{N_1}$-periodic function of the form
$$
\frac{\pi n \mathrm{i}}{N_1}\,\big( 1 -  {\mathrm e}^{-2\pi {\mathrm i}\, n/N_1}\big)^{-1}\,  {\mathrm e}^{-2\pi {\mathrm i}\, n\,x}\,, \quad x \in (0,\,\frac{1}{N_1})\,.
$$
\end{Lemma}

\emph{Proof.} For fixed $n\in I_N \setminus \{0\}$, we consider the Fourier series of the special $\frac{1}{N_1}$-periodic function
$g_n(x) := \mathrm{e}^{-2 \pi \mathrm{i}\,nx}$ for $x \in (0,\,\frac{1}{N_1})$
with
$$
g_n(0) = g_n(\frac{1}{N_1}) := \frac{1}{2}\,\big( 1 +  \mathrm{e}^{-2 \pi \mathrm{i}\,n/N_1}\big)\,.
$$
For $n=0$ we have $g_0(x) = 1$. Then the $r$th Fourier coefficient of $g_n$ reads as follows
$$
N_1\,\int_0^{1/N_1} g_n(t)\,\mathrm{e}^{-2 \pi \mathrm{i}\,r N_1\,t}\,\mathrm{d}t = \frac{1}{\pi \mathrm{i}}\,\big( 1 -  \mathrm{e}^{-2 \pi \mathrm{i}\,n/N_1}\big)\,\frac{1}{r + \frac{n}{N_1}}
$$
for $r \in \mathbb Z$.
By the Convergence Theorem of Dirichlet--Jordan, the Fourier series of $g_n$ is pointwise convergent such that for each $x \in \mathbb R$ and $n\in I_N \setminus \{0\}$
\begin{equation*} g_n(x) = \frac{1}{\pi \mathrm{i}}\,\big( 1 -  \mathrm{e}^{-2 \pi \mathrm{i}\,n/N_1}\big)\,\sum_{r\in \mathbb Z}\frac{1}{r + \frac{n}{N_1}}\,\mathrm{e}^{2 \pi \mathrm{i}\,r N_1 \,x}\,.
\end{equation*}
This completes the proof. \qedsymbol

\begin{Lemma}
\label{Lemma:errorrectwindow}
The $C(\mathbb T)$-error constant of the rectangular window function \eqref{eq:rectwindow} can be estimated by
\begin{equation}
\label{eq:errorrectwindow}
0.18 < \frac{1}{2}-\frac{1}{\pi} \le e_{\sigma}({\varphi}_{\mathrm{rect}}) \le \frac{1}{2} + \frac{\pi}{4} < 1.3\, ,
\end{equation}
i.e., the rectangular window function \eqref{eq:rectwindow}  is not convenient for $\mathrm{NFFT}$.
\end{Lemma}

\emph{Proof.} By \eqref{eq:FTrectwindow} we obtain for $n=0$ and $r \in \mathbb Z$ that
$$
\frac{{\hat\varphi}_{\mathrm{rect}}(r N_1)}{{\hat\varphi}_{\mathrm{rect}}(0)} = \left\{\begin{array}{ll} 1 & \quad r = 0\,,\\
0 & \quad r \in \mathbb Z \setminus \{0\}
\end{array}\right.
$$
and hence
$$
\sum_{r\in \mathbb Z\setminus \{0\}} \frac{{\hat\varphi}_{\mathrm{rect}}(r N_1)}{{\hat\varphi}_{\mathrm{rect}}(0)}\,{\mathrm e}^{2\pi {\mathrm i}\, r N_1\,x} = 0\,.
$$
By \eqref{eq:FTrectwindow} we have for $n \in I_N \setminus \{0\}$ and $r \in \mathbb Z$,
$
\frac{{\hat\varphi}_{\mathrm{rect}}(n + r N_1)}{{\hat\varphi}_{\mathrm{rect}}(n)} = \frac{n}{n + r N_1}\,.
$
Thus we obtain by Lemma \ref{Lemma:FTrectwindow} that for $x\in (0,\,\frac{1}{N_1})$ and $n \in I_N \setminus \{0\}$
\begin{eqnarray*}
\sum_{r\in \mathbb Z\setminus \{0\}} \frac{{\hat\varphi}_{\mathrm{rect}}(n + r N_1)}{{\hat\varphi}_{\mathrm{rect}}(n)}
\,{\mathrm e}^{2\pi {\mathrm i}\, r N_1\,x}
&=&  \frac{n}{N_1}\, \sum_{r\in \mathbb Z\setminus \{0\}} \frac{1}{r + \frac{n}{N_1}}\,{\mathrm e}^{2\pi {\mathrm i}\, r N_1\,x}\\
&=&  \frac{\pi n \mathrm{i}}{N_1}\,\big(1 -  \mathrm{e}^{-2 \pi \mathrm{i}\,n/N_1}\big)^{-1}\,{\mathrm e}^{-2\pi {\mathrm i}\,n x} - \frac{n}{N_1}\,.
\end{eqnarray*}
Using
$
\big| 1 - \mathrm{e}^{-2 \pi \mathrm{i}\,n/N_1}\big| = 2 \,\big|\sin \frac{\pi n}{N_1}\big| \,,
$
it follows that
\begin{eqnarray*}
& &\big| \frac{\pi n \mathrm{i}}{N_1}\,\big(1 -  \mathrm{e}^{-2 \pi \mathrm{i}\,n/N_1}\big)^{-1}\,{\mathrm e}^{-2\pi {\mathrm i}\,n x} - \frac{n}{N_1}\big| \ge \frac{|n|}{N_1}\,\big(\frac{\pi}{2}\,|\sin\frac{\pi n}{N_1}|^{-1} - 1\big)\\
& & = \frac{|n|}{2 N_1}\,\frac{\pi - 2\,|\sin\frac{\pi n}{N_1}|}{|\sin\frac{\pi n}{N_1}| }
\ge \frac{1}{2\pi}\,\big(\pi - 2\,|\sin \frac{\pi n}{\sigma N}|\big) \ge \frac{1}{2\pi}\,(\pi - 2)\,.
\end{eqnarray*}
Analogously, we estimate
\begin{eqnarray*}
& &\big| \frac{\pi n \mathrm{i}}{N_1}\,\big(1 -  \mathrm{e}^{-2 \pi \mathrm{i}\,n/N_1}\big)^{-1}\,{\mathrm e}^{-2\pi {\mathrm i}\,n x} - \frac{n}{N_1}\big| \le \frac{|n|}{N_1}\,\big(\frac{\pi}{2}\,|\sin\frac{\pi n}{N_1}|^{-1} + 1\big)\\
& & = \frac{|n|}{2 N_1}\,\frac{\pi + 2\,|\sin\frac{\pi n}{N_1}|}{|\sin\frac{\pi n}{N_1}| }
\le \frac{1}{4}\,\big(\pi + 2\,|\sin \frac{\pi n}{N_1}|\big) \le \frac{1}{4}\,(\pi + 2)\,.
\end{eqnarray*}
Consequently the rectangular window function is not convenient for NFFT, since the corresponding $C(\mathbb T)$-error constant $e_{\sigma}(\varphi_{\mathrm{rect}})$ can be estimated by \eqref{eq:errorrectwindow}. \qedsymbol

\begin{Lemma}
\label{Lemma:sumFTvarphirect}
The Fourier transform of rectangular window function \eqref{eq:rectwindow} has the following property
\begin{equation}
\label{eq:sumhatphirect1}
\big(1 - \frac{2}{\pi}\big)\,\frac{m}{N_1}\,\big|{\mathrm{sinc}}\frac{\pi m}{\sigma}\big| \le \big\|\sum_{r \in \mathbb Z \setminus \{0\}} {\hat \varphi}_{\mathrm{rect}}(n + r N_1)\,{\mathrm e}^{2\pi {\mathrm i}r N_1 \, \cdot}\big\|_{C(\mathbb T)} \le \frac{3m}{N_1}\,.
\end{equation}
\end{Lemma}

\emph{Proof}. By \eqref{eq:convenient1} -- \eqref{eq:esigmaNvarphi}, we obtain
$$
\max_{n \in I_N} \big\|\sum_{r \in \mathbb Z \setminus \{0\}} {\hat \varphi}_{\mathrm{rect}}(n + r N_1)\,{\mathrm e}^{2\pi {\mathrm i}r N_1 \, \cdot}\big\|_{C(\mathbb T)} \ge e_{\sigma}( {\varphi}_{\mathrm{rect}})\,
\min_{n \in I_N} | {\hat \varphi}_{\mathrm{rect}}(n)|\,.
$$
Then from Lemma \ref{Lemma:errorrectwindow} and \eqref{eq:FTrectwindow} it follows immediately the lower estimate in \eqref{eq:sumhatphirect1}.

Now we show the upper estimate in \eqref{eq:sumhatphirect1}.
For each $n \in I_N \setminus \{0\}$, the $\frac{1}{N_1}$-periodic Fourier series
\begin{eqnarray*}
f_n(x)&:=& \sum_{r \in \mathbb Z} {\hat \varphi}_{\mathrm{rect}}(n + r N_1)\,{\mathrm e}^{2\pi {\mathrm i}r N_1 x} = \frac{2m}{N_1}\,\sum_{r \in \mathbb Z} \mathrm{sinc}\,\frac{2\pi m\,(n + r N_1)}{N_1}\,{\mathrm e}^{2\pi {\mathrm i}r N_1 x}\\
&=& \frac{1}{\pi N_1}\,\sin \frac{2 \pi m n}{N_1}\, \sum_{r \in \mathbb Z} \frac{1}{r + \frac{n}{N_1}}\,{\mathrm e}^{2\pi {\mathrm i}r N_1 x}
\end{eqnarray*}
is pointwise convergent by Lemma \ref{Lemma:FTrectwindow} such that
$$
f_n(x) = \left\{ \begin{array}{ll}
\frac{\mathrm i}{2m}\,\sin\frac{2 \pi m n}{N_1}\,\big(1 - {\mathrm e}^{- 2 \pi {\mathrm i} n /N_1}\big)^{-1}\, {\mathrm e}^{-2\pi {\mathrm i}n x} &\; x \in \big(0,\,\frac{1}{N_1}\big)\,,\\ [1ex]
\frac{1}{4m}\,\sin\frac{2 \pi m n}{N_1}\,\cot \frac{\pi n}{N_1} & \; x \in \{0,\,\frac{1}{N_1}\}\,.
\end{array} \right.
$$
In the case $n=0$ we have
$
f_0(x) := \sum_{r \in \mathbb Z} \mathrm{sinc} (2\pi m  r)\,{\mathrm e}^{2\pi {\mathrm i}r N_1 x} = 1\, .
$
For arbitrary $x \in \mathbb R$ and $m \in \mathbb N$, it holds obviously
$
\big| \sin (2 m x)\big| \le 2 m \,|\sin x|
$
and so
$
\big| 1 -  {\mathrm e}^{- 2 \pi {\mathrm i} n /N_1}\big| = 2 \, \big| \sin \frac{\pi n}{N_1}\big|\, .
$
We obtain for $n \in I_N \setminus \{0\}$ that
$$
\max_{x \in \mathbb R} |f_n(x)| \le \frac{1}{2 N_1}\,\big| \sin \frac{2\pi m n}{N_1}\big|\,\big| \sin \frac{\pi n}{N_1}\big|^{-1} \le \frac{m}{N_1}\,.
$$
Thus for $n \in I_N \setminus \{0\}$ and $x\in \big[0,\,\frac{1}{N_1}\big]$, we receive
$$
\sum_{r \in \mathbb Z \setminus \{0\}} {\hat \varphi}_{\mathrm{rect}}(n + r N_1)\,{\mathrm e}^{2\pi {\mathrm i}r N_1 x} = f_n(x) - \frac{2 m}{N_1}\,\mathrm{sinc} \frac{2 \pi m n}{N_1}
$$
and hence
$$
\big\|\sum_{r \in \mathbb Z \setminus \{0\}} {\hat \varphi}_{\mathrm{rect}}(n + r N_1)\,{\mathrm e}^{2\pi {\mathrm i}r N_1 \, \cdot}\big\|_{C(\mathbb T)} \le \frac{3m}{N_1}\,.
$$
In the case $n=0$, the above estimate is also true, since
$$
\sum_{r \in \mathbb Z \setminus \{0\}} {\hat \varphi}_{\mathrm{rect}}(r N_1)\,{\mathrm e}^{2\pi {\mathrm i}r N_1 x} = 0\,.
$$
This completes the proof. \qedsymbol

\section{Continuous Kaiser--Bessel window function}\label{Sec:KBWindow}

In the following, we consider the so-called \emph{continuous Kaiser--Bessel window function}
\begin{equation}
\label{eq:Kaiser-Bessel}
\varphi_{\mathrm{cKB}}(x) := \left\{ \begin{array}{ll}  \frac{1}{I_0(\beta)-1}\Big(I_0\big(\beta\, \sqrt{1 - \big(\frac{N_1 x}{m}\big)^2}\big) -1\Big) & \quad x \in I\,,\\ [1ex]
0 & \quad x \in \mathbb R \setminus  I
\end{array}\right.
\end{equation}
with $b:= 2\pi (1 - \frac{1}{2\sigma})$, $\sigma > 1$, $\beta := bm$, $N\in 2 \,\mathbb N$, and $N_1 \in 2 \,\mathbb N$, see \cite{Ne14}, where
$$
I_0(x) := \sum_{k=0}^{\infty} \frac{1}{(k!)^2}\,\big(\frac{x}{2}\big)^{2k}\,, \quad x \in \mathbb R\,,
$$
denotes the \emph{modified Bessel function of first kind}. Further we assume that $m \in \mathbb N \setminus \{1\}$ fulfills  $2 m \ll N_1$.
We emphasize that the \emph{shape parameter} $\beta$ in \eqref{eq:Kaiser-Bessel} depends on $m$ and $\sigma$.
Obviously, $\varphi_{\mathrm{cKB}} \in \Phi_{m,N_1}$ is a continuous window function. Note that
a discontinuous version \eqref{eq:dKaiser-Bessel} of \eqref{eq:Kaiser-Bessel} (see Remark \ref{Remark:KB-window}) was considered in \cite{Fou02}, but the $C(\mathbb T)$-error constant of \eqref{eq:dKaiser-Bessel} was not determined.

By \cite[p.~95]{Oberh90} and \eqref{eq:FTrectwindow}, the Fourier transform of the continuous Kaiser--Bessel window function \eqref{eq:Kaiser-Bessel} has the form
\begin{equation*} {\hat \varphi}_{\mathrm{cKB}}(v) =   \frac{2m}{(I_0(\beta)-1)\,N_1}\,\Big(\frac{\sinh \big(\beta\,\sqrt{1- \big(\frac{2\pi v}{N_1 b}\big)^2}\big)}{\beta \,\sqrt{1- \big(\frac{2\pi v}{N_1 b}\big)^2}} - {\mathrm{sinc}}\frac{2 \pi m v}{N_1} \Big)
\end{equation*}
for $|v| < \frac{N_1 b}{2 \pi}$ and
\begin{equation*} {\hat \varphi}_{\mathrm{cKB}}(v) =   \frac{2m}{(I_0(\beta)-1)\,N_1}\,\Big({\mathrm{sinc}}\big(\beta\,\sqrt{\big(\frac{2\pi v}{N_1 b}\big)^2-1}\big) - {\mathrm{sinc}}\frac{2 \pi m v}{N_1} \Big)
\end{equation*}
for $|v| \ge \frac{N_1 b}{2\pi}$.
One can show that ${\hat \varphi}_{\mathrm{cKB}} |_{[0,\,\frac{N_1 b}{2\pi})}$ is positive and decreasing such that
\begin{equation}
\label{eq:phiKB(N/2)}
\min_{n \in I_N}{\hat \varphi}_{\mathrm{cKB}}(n) = {\hat \varphi}_{\mathrm{cKB}}\big(\frac{N}{2}\big) \ge \frac{2 m}{(I_0(\beta) - 1)\,N_1}\Big[ \frac{\sinh \big(2\pi m\sqrt{1 - 1/\sigma}\big)}{2m \pi \sqrt{1 - 1/\sigma}} - \frac{\sigma}{\pi m}\Big]\,.
\end{equation}
Using the scaled frequency $w = 2\pi m v/N_1$, we obtain
\begin{equation}
\label{eq:FTKaiser-Besselinw}
{\hat \varphi}_{\mathrm{cKB}}\big(\frac{N_1 w}{2\pi m}\big) = \frac{2m}{(I_0(\beta)-1)\,N_1} \cdot \left\{ \begin{array}{ll} \Big(\frac{\sinh (\sqrt{\beta^2- w^2})}{\sqrt{\beta^2- w^2}}- {\mathrm{sinc}}\,w\Big) & \quad |w| < \beta\,,\\ [2ex]
\big( {\mathrm{sinc}}(\sqrt{w^2-\beta^2}) - {\mathrm{sinc}}\,w\big) & \quad |w| \ge \beta\,.
\end{array}\right.
\end{equation}
Note that we have
$
\frac{N_1 b}{2\pi} = N_1 - \frac{N}{2} > \frac{N}{2}
$
by the special choice of $b$.

\begin{Lemma}
\label{Lemma:sincsqrt-sinc}
For $|w|\ge \beta$ we have
$$
\big|{\mathrm{sinc}}(\sqrt{w^2-\beta^2}) - {\mathrm{sinc}}\,w\big| \le \frac{2 \beta^2}{w^2}\,.
$$
\end{Lemma}

\emph{Proof}. Since the $\mathrm{sinc}$-function is even, we consider only the case $w\ge \beta$. For $w = \beta$, the above inequality is true, since
$
|{\mathrm{sinc}}\,0 - {\mathrm{sinc}}\,\beta| \le 1 + |{\mathrm{sinc}}\,\beta| < 2\,.
$
For $w > \beta$ we obtain
\begin{eqnarray*}
& &\big|{\mathrm{sinc}}(\sqrt{w^2-\beta^2}) - {\mathrm{sinc}}\,w\big| = \big| \frac{\sin\sqrt{w^2 - \beta^2}}{\sqrt{w^2 - \beta^2}} - \frac{\sin w}{w}\big|\\
& & \le \frac{1}{w}\,\big| \sin \sqrt{w^2 - \beta^2} - \sin w \big| + \big| \sin \sqrt{w^2 - \beta^2}\big| \big(\frac{1}{\sqrt{w^2- \beta^2}} - \frac{1}{w}\big)\\
& & \le \frac{2}{w}\,\big|\sin \frac{1}{2}\big(\sqrt{w^2 - \beta^2} - w\big)\big| + \big| \sin \sqrt{w^2 - \beta^2}\big| \big(\frac{1}{\sqrt{w^2-\beta^2}} - \frac{1}{w}\big)\,.
\end{eqnarray*}
From
\begin{equation}
\label{eq:w-sqrt}
\sqrt{w^2 - \beta^2} - w = w\,\big(1 - \sqrt{1 - \frac{\beta^2}{w^2}}\big)= w\,\frac{1 - \big(1- \frac{\beta^2}{w^2}\big)}{1 + \sqrt{1-\frac{\beta^2}{w^2}}}\le \frac{\beta^2}{w}
\end{equation}
it follows that
$$
\big|\sin \frac{1}{2}\big(\sqrt{w^2 - \beta^2} - w\big)\big| \le  \frac{1}{2}\big(\sqrt{w^2 - \beta^2} - w\big) \le \frac{\beta^2}{2w}\,.
$$
Further we receive by \eqref{eq:w-sqrt} that
\begin{eqnarray*}
& &\big| \sin \sqrt{w^2 - \beta^2}\big| \big(\frac{1}{\sqrt{w^2- \beta^2}} - \frac{1}{w}\big) = \big| \sin \sqrt{w^2 - \beta^2}\big|\,\frac{w -\sqrt{w^2- \beta^2}}{w\,\sqrt{w^2 - \beta^2}} \\
& &  \le \big| \sin \sqrt{w^2 - \beta^2}\big|\,\frac{\beta^2}{w^2} \le \frac{\beta^2}{w^2}\,.
\end{eqnarray*}
This completes the proof. \qedsymbol
\medskip

In our study we use the following

\begin{Lemma}
\label{Lemma:sumf(u+r)}
Let $f:\,(0,\,\infty) \to (0,\,\infty)$ be a decreasing function which is integrable on each interval $[1 - |u|,\,\infty)$ with arbitrary $u \in (-1,\,1)$. Further we extend $f$ by $f(-x) := f(x)$ for all $x > 0$.

Then for each $u \in (-1,\,1)$, it follows
\begin{eqnarray*}
\sum_{r\in {\mathbb Z}\setminus \{0,\pm 1\}} f(u+r) &\le& 2 \,\int_{1 - |u|}^{\infty} f(t)\,{\mathrm d}t\,,\\
\sum_{r\in {\mathbb Z}\setminus \{0\}} f(u+r) &\le& 2 \,f(1 - |u|) + 2 \,\int_{1 - |u|}^{\infty} f(t)\,{\mathrm d}t\,.
\end{eqnarray*}
\end{Lemma}

\emph{Proof}. For arbitrary $u \in (-1,\,1)$ and $r \in \mathbb N$, we have
\begin{equation}
\label{eq:f(u+r)}
f(u + r) \le f(r - |u|)\,.
\end{equation}
Using \eqref{eq:f(u+r)}, the following series can be estimated by
\begin{eqnarray*}
\sum_{r=2}^{\infty} f(u + r) &\le& \sum_{r=2}^{\infty} f(r - |u|)\,,\\
\sum_{r=2}^{\infty} f(u - r) &\le& \sum_{r=2}^{\infty} f(r - |u|)\,.
\end{eqnarray*}
Hence it follows by the integral test for convergence of series that
\begin{eqnarray*}
\sum_{r\in {\mathbb Z}\setminus \{0,\pm 1\}} f(u+r) &=& \sum_{r=2}^{\infty} f(u + r) + \sum_{r=2}^{\infty} f(u - r)\le 2\,\sum_{r=2}^{\infty} f(r - |u|)\\
&\le& 2\,\int_1^{\infty} f(x - |u|)\,{\mathrm d}x = 2 \,\int_{1 - |u|}^{\infty} f(t)\,{\mathrm d}t\,.
\end{eqnarray*}
Hence we conclude that
$$
\sum_{r\in {\mathbb Z}\setminus \{0\}} f(u+r) \le
 2 \,f(1 - |u|) + 2 \,\int_{1 - |u|}^{\infty} f(t)\,{\mathrm d}t\,. \quad \qedsymbol
$$

We illustrate Lemma \ref{Lemma:sumf(u+r)} for some special functions $f$, which we need later.

\begin{Example}\label{Example:specialsumf(u+r)}
For the function $f(x) = x^{-\mu}$, $x>0$, with $\mu > 1$, Lemma \ref{Lemma:sumf(u+r)} provides that for each $u \in (-1,\,1)$ it holds
\begin{equation}
\label{eq:sumxmu}
\sum_{r\in {\mathbb Z}\setminus \{0\}} |u+r|^{-\mu} \le 2 \,(1 - |u|)^{-\mu} + \frac{2\,(1 - |u|)^{1- \mu}}{\mu - 1}\,.
\end{equation}
Especially for $\mu = 2$, it follows that
\begin{equation}
\label{eq:sumx2}
\sum_{r\in {\mathbb Z}\setminus \{0\}} |u+r|^{-2} \le \frac{4 - 2\,|u|}{(1 - |u|)^2}\le \frac{4}{(1 - |u|)^2}\,.
\end{equation}
For the function $f(x) = {\mathrm e}^{-a x}$, $x>0$, with $a > 0$, we obtain by Lemma \ref{Lemma:sumf(u+r)} that for each $u \in (-1,\,1)$,
\begin{equation}
\label{eq:sumexp(-ax)}
\sum_{r\in {\mathbb Z}\setminus \{0\}} {\mathrm e}^{- a\,|u + r|} \le \big(2 + \frac{2}{a}\big)\,{\mathrm e}^{a\,|u| -a}\,.
\end{equation}
Choosing the function $f(x) = \frac{1}{ax}\,{\mathrm e}^{- \sqrt{a\,x}}$, $x > 0$, with $a > 0$, Lemma \ref{Lemma:sumf(u+r)} implies that for each $u\in (-1,\,1)$,
\begin{equation}
\label{eq:sumexp(-sqrtx)/x}
\sum_{r\in {\mathbb Z}\setminus \{0\}} \frac{1}{a\,|u+r|}\,{\mathrm e}^{- \sqrt{ a\,|u + r|}} \le \frac{2}{a\,(1 -|u|)}\,{\mathrm e}^{-\sqrt{a - a\,|u|}} + \frac{4}{a}\,E_1\big(\sqrt{a\,(1-|u|)}\big)\,,
\end{equation}
where $E_1$ denotes the \emph{exponential integral}
$
E_1(x) := \int_x^{\infty} \frac{1}{t}\,{\mathrm e}^{-t}\,{\mathrm d}t$, $x > 0$.
\end{Example}

\begin{Lemma}
\label{Lemma4.2}
For all $n \in I_N$ we estimate
$$
\sum_{r \in \mathbb Z \setminus \{0\}} \big| {\hat \varphi}_{\mathrm{cKB}}(n + r N_1)\big| \le \frac{8m}{(I_0(\beta)-1)\,N_1}\,.
$$
\end{Lemma}

\emph{Proof}. By \eqref{eq:FTKaiser-Besselinw} and Lemma \ref{Lemma:sincsqrt-sinc} we obtain for all frequencies $|v| \ge \frac{N_1 b}{2 \pi} = N_1 - \frac{N}{2}$ that
$$
\big| {\hat \varphi}_{\mathrm{cKB}}(v) \big| \le \frac{m b^2 N_1}{2\,(I_0(bm)-1)\, \pi^2 v^2}\,.
$$
Thus we have for each $n\in I_N$
\begin{eqnarray*}
& &\sum_{r \in \mathbb Z \setminus \{0\}} \big| {\hat \varphi}_{\mathrm{cKB}}(n + r N_1)\big| \le \frac{b \beta}{2\,(I_0(\beta)- 1)\, N_1 \pi^2} \,\sum_{r \in \mathbb Z \setminus \{0\}} \big(r + \frac{n}{N_1}\big)^{-2}\\
&  &\le \, \frac{2 b \beta}{(I_0(\beta) - 1)\,N_1 \pi^2}\, \big(1 - \frac{|n|}{N_1}\big)^{-2}\le \frac{2 b \beta}{(I_0(\beta) - 1)\,N_1 \pi^2}\, \big(1 - \frac{1}{2\sigma}\big)^{-2}\,,
\end{eqnarray*}
since by \eqref{eq:sumx2} it holds for all $n \in I_N$,
$$
\sum_{r \in \mathbb Z \setminus \{0\}} \big(r + \frac{n}{N_1}\big)^{-2} \le 4\,\big(1 - \frac{|n|}{N_1}\big)^{-2} \le  4\,\big(1 - \frac{1}{2\sigma}\big)^{-2}\,.
$$
By the special choice of $b = 2\pi \,\big(1 - \frac{1}{2\sigma}\big)$, we obtain the above inequality. \qedsymbol
\medskip

From Lemma \ref{Lemma4.2} it follows immediately that
for all $n \in I_N$,
\begin{equation}
\label{eq:sumphiKB}
\big\|\sum_{r \in \mathbb Z \setminus \{0\}}  {\hat \varphi}_{\mathrm{cKB}}(n + r N_1)\,{\mathrm e}^{2 \pi {\mathrm i} r N_1 \, \cdot}\big\|_{C(\mathbb T)} \le \sum_{r \in \mathbb Z \setminus \{0\}} \big| {\hat \varphi}_{\mathrm{cKB}}(n + r N_1)\big| \le \frac{8 m}{(I_0(\beta) - 1)\,N_1}\,.
\end{equation}

\begin{Theorem}
\label{Theorem4.4}
Let $b = 2\pi\,(1 - \frac{1}{2\sigma})$, $\sigma > 1$, $\beta = b m$, $N \in 2 \mathbb N$, and $N_1 = \sigma N\in 2 \mathbb N$. Further $m\in \mathbb N \setminus \{1\}$ with $2m \ll N_1$ is  given.

Then the $C(\mathbb T)$-error constant of the continuous Kaiser--Bessel window function \eqref{eq:Kaiser-Bessel} can be estimated by
\begin{equation}
\label{eq:esigmaKB}
e_{\sigma}(\varphi_{\mathrm{cKB}}) \le 16\, m \pi \, \sqrt{1-1/\sigma}\,\Big[{\mathrm e}^{2\pi m \sqrt{1- 1/\sigma}} - {\mathrm e}^{-2 \pi m\,\sqrt{1 - 1/\sigma}}
- 4 \sqrt{\sigma^2 - \sigma}
\Big]^{-1}\,.
\end{equation}
\end{Theorem}

Note that for $\sigma \ge \frac{5}{4}$, it holds $2 \pi\,\sqrt{1 - 1/\sigma} \ge 2\pi/\sqrt 5$ and hence
\begin{equation}
\label{eq:sigma5/4}
{\mathrm e}^{-2 \pi m\,\sqrt{1 - 1/\sigma}} \le {\mathrm e}^{-2\pi m/\sqrt 5} < 0.06^m\,.
\end{equation}
\emph{Proof.} By the definition \eqref{eq:convenient1} of the $C(\mathbb T)$-error constant, it holds
$$
e_{\sigma}(\varphi_{\mathrm{cKB}}) = \sup_{N \in 2 \mathbb N} e_{\sigma, N}(\varphi_{\mathrm{cKB}})
$$
with
\begin{eqnarray*}
e_{\sigma, N}(\varphi_{\mathrm{cKB}}) &=& \max_{n \in I_N} \big\|\sum_{r \in \mathbb Z \setminus \{0\}}  \frac{{\hat \varphi}_{\mathrm{cKB}}(n + r N_1)}{{\hat \varphi}_{\mathrm{cKB}}(n)}\,{\mathrm e}^{2 \pi {\mathrm i} r N_1 \, \cdot}
\big\|_{C(\mathbb T)}\\ [1ex]
&\le& \frac{1}{\min_{n \in I_N}{\hat \varphi}_{\mathrm{cKB}}(n)}\, \max_{n \in I_N} \big\|\sum_{r \in \mathbb Z \setminus \{0\}}  {\hat \varphi}_{\mathrm{cKB}}(n + r N_1)\,{\mathrm e}^{2 \pi {\mathrm i} r N_1 \,\cdot}
\big\|_{C(\mathbb T)}\,,
\end{eqnarray*}
where it holds \eqref{eq:phiKB(N/2)}, i.e.,
\begin{eqnarray*}
{\hat \varphi}_{\mathrm{cKB}}\big(\frac{N}{2}\big) &\ge& \frac{1}{(I_0(\beta) - 1)\,N_1 \pi\,\sqrt{1-1/\sigma}}\,\big[\sinh\big(2\pi m \sqrt{1 - 1/\sigma}\big) - 2\, \sqrt{\sigma^2 - \sigma}\big]\\
&=& \frac{1}{2\,(I_0(\beta) - 1)\,N_1 \pi\,\sqrt{1-1/\sigma}}\,\big[{\mathrm e}^{2\pi m \sqrt{1 - 1/\sigma}} - {\mathrm e}^{-2\pi m \sqrt{1 - 1/\sigma}} - 4\, \sqrt{\sigma^2 - \sigma}\big]\,.
\end{eqnarray*}
Thus from \eqref{eq:sumphiKB} it follows the assertion \eqref{eq:esigmaKB}. \qedsymbol
\medskip

\begin{Remark}
\label{Remark:KB-window}
As in \cite{Fou02} and \cite{Ba20}, we consider also the \emph{standard Kaiser--Bessel window function}
\begin{equation}
\label{eq:dKaiser-Bessel}
\varphi_{\mathrm{KB}}(x) := \left\{ \begin{array}{ll}  \frac{1}{I_0(\beta)}\, I_0\big(\beta\, \sqrt{1 - \big(\frac{N_1 x}{m}\big)^2}\big) & \quad x \in I\,,\\ [1ex]
\frac{1}{2\,I_0(\beta)} & \quad x = \pm \frac{m}{N_1}\,, \\ [1ex]
0 & \quad x \in \mathbb R \setminus  \bar I
\end{array}\right.
\end{equation}
with the shape parameter $\beta = m b = 2\pi m\, (1 - \frac{1}{2\sigma})$, $\sigma > 1$, $N\in 2 \,\mathbb N$, and $N_1 \in 2 \,\mathbb N$. Further we assume that $m \in \mathbb N \setminus \{1\}$
fulfills $2 m \ll N_1$. This window function possesses jump discontinuities at $x = \pm \frac{m}{N_1}$ with very small jump height $I_0(\beta)^{-1}$, such that \eqref{eq:dKaiser-Bessel}
is ``almost continuous''. The Fourier transform of \eqref{eq:dKaiser-Bessel} is even and reads by \cite[p.~95]{Oberh90} as follows
\begin{equation*}
{\hat \varphi}_{\mathrm{KB}}(v) =  \left\{ \begin{array}{ll}
\frac{2m}{I_0(\beta)\,N_1}\,\frac{\sinh \big(\beta\,\sqrt{1- \big(\frac{2\pi v}{N_1 b}\big)^2}\big)}{\beta \,\sqrt{1- \big(\frac{2\pi v}{N_1 b}\big)^2}} & \quad |v| < \frac{N_1 b}{2 \pi}\,, \\ [1ex]
 \frac{2m}{I_0(\beta)\,N_1}\,{\mathrm{sinc}}\big(\beta\,\sqrt{\big(\frac{2\pi v}{N_1 b}\big)^2-1}\big) & \quad |v| \ge \frac{N_1 b}{2\pi}\,.
\end{array} \right.
\end{equation*}
Thus ${\hat \varphi}_{\mathrm{KB}}|_{[0,\,\frac{N_1 b}{2\pi})}$ is positive and decreasing such that
$$
\min_{n \in I_N} {\hat \varphi}_{\mathrm{KB}}(n) = {\hat \varphi}_{\mathrm{KB}}\big(\frac{N}{2}\big) \ge \frac{1}{I_0(\beta)\,N_1 \pi}\,\frac{\sinh \big(2 \pi m\,\sqrt{1- 1/\sigma}\big)}{\sqrt{1- 1/\sigma}}\,.
$$
Splitting ${\varphi}_{\mathrm{KB}}$ in the form
$$
{\varphi}_{\mathrm{KB}}(x) = \big(1 - \frac{1}{I_0(\beta)}\big)\,{\varphi}_{\mathrm{cKB}}(x) + \frac{1}{I_0(\beta)}\,{\varphi}_{\mathrm{rect}}(x)\,, \quad x \in \mathbb R\,,
$$
we can estimate for all $n \in I_N$,
\begin{eqnarray*}
\big\|\sum_{r \in \mathbb Z \setminus \{0\}} {\hat \varphi}_{\mathrm{KB}}(n + r N_1)\,{\mathrm e}^{2\pi {\mathrm i}r N_1 \, \cdot}\big\|_{C(\mathbb T)}
&\le&  \big(1 - \frac{1}{I_0(\beta)}\big)\,\sum_{r \in \mathbb Z \setminus \{0\}} \big|{\hat \varphi}_{\mathrm{cKB}}(n + r N_1)\big|\\
& & +\, \frac{1}{I_0(\beta)}\, \big\|\sum_{r \in \mathbb Z \setminus \{0\}} {\hat \varphi}_{\mathrm{rect}}(n + r N_1)\,{\mathrm e}^{2\pi {\mathrm i}r N_1 \,\cdot}\big\|_{C(\mathbb T)}\,.
\end{eqnarray*}
Then from Lemma \ref{Lemma4.2} and \eqref{eq:sumhatphirect1} it follows that
$$
\max_{n \in I_N} \big\|\sum_{r \in \mathbb Z \setminus \{0\}} {\hat \varphi}_{\mathrm{KB}}(n + r N_1)\,{\mathrm e}^{2\pi {\mathrm i}r N_1 \, \cdot}\big\|_{C(\mathbb T)} \le \frac{11 m}{I_0(\beta)\,N_1}\,.
$$
Consequently, we obtain the following estimate of the $C(\mathbb T)$-error constant of \eqref{eq:dKaiser-Bessel}
\begin{eqnarray*}
e_{\sigma}(\varphi_{\mathrm{KB}}) &\le& 11\, m \pi \, \sqrt{1-\frac{1}{\sigma}}\,\big[\sinh(2\pi m \sqrt{1- 1/\sigma})\big]^{-1} \nonumber \\ [1ex]
&=& 22 \pi m \, \sqrt{1-\frac{1}{\sigma}}\,\Big[{\mathrm e}^{2\pi m \sqrt{1- 1/\sigma}} - {\mathrm e}^{-2 \pi m \, \sqrt{1 -1/\sigma}} \Big]^{-1}\,.
\end{eqnarray*}
For $\sigma \ge \frac{5}{4}$, we sustain by \eqref{eq:sigma5/4} that
$$
e_{\sigma}(\varphi_{\mathrm{KB}}) \le 22 \pi m \, \sqrt{1-\frac{1}{\sigma}}\,\Big[{\mathrm e}^{2\pi m \sqrt{1- 1/\sigma}} - 0.06^m \Big]^{-1}\,.
$$
\end{Remark}

\section{sinh-type window function}
\label{Sec:sinhWindow}

For fixed shape parameter $\beta = b m$ with $m\in \mathbb N \setminus \{1\}$, $b = 2 \pi\,\big(1 - \frac{1}{2 \sigma}\big)$, and oversampling factor $\sigma \ge \frac{5}{4}$, we consider the $\sinh$-\emph{type window function}
\begin{equation}
\label{eq:sinhwindow}
\varphi_{\sinh}(x) := \left\{ \begin{array}{ll}
\frac{1}{\sinh \beta}\,\sinh \big({\beta \,\sqrt{1 - (N_1x/m)^2}}\big) & \quad x \in I\,, \\ [1ex]
0 & \quad x \in \mathbb R \setminus  I\,.
\end{array}\right.
\end{equation}
Obviously, $\varphi_{\sinh} \in  \Phi_{m,N_1}$ is a continuous window function which was introduced by the authors in \cite{PoTa20} (with another shape parameter $\beta = 4 m$). Let $N \ge 8$ be an even integer.
For a sampling factor $\sigma \ge \frac{5}{4}$
and $b = 2\pi \,(1 - \frac{1}{2 \sigma})$, we form $N_1 = \sigma N \in 2 \mathbb N$ and the shape parameter $\beta = b m$, where $m \in \mathbb N \setminus \{1\}$ with $2m \ll N_1$.

Substituting $t = N_1 x /m$, we determine the even Fourier transform
$$
{\hat \varphi}_{\sinh}(v) = \int_{\mathbb R} \varphi_{\sinh}(x)\,{\mathrm e}^{-2\pi {\mathrm i}\,v x}\,{\mathrm d}x
= \frac{2m}{N_1\,\sinh \beta}\,\int_0^1 \sinh \big(\beta\,\sqrt{1 - t^2}\big)\,\cos\frac{2\pi\,m v t}{N_1}\,{\mathrm d}t\,.
$$
Using the scaled frequency
$w := 2 \pi m v/N_1$, the Fourier transform of \eqref{eq:sinhwindow} reads by \cite[p.~38]{Oberh90} as follows
\begin{equation}
\label{eq:FTsinhwindow}
{\hat \varphi}_{\sinh}\big(\frac{N_1 w}{2\pi m}\big) = \frac{\pi m \beta}{N_1\, \sinh \beta}\,\left\{ \begin{array}{ll}
(\beta^2 - w^2)^{-1/2}\,I_1\big(\sqrt{\beta^2 - w^2}\big) &\quad w\in (- \beta,\,\beta)\,,\\
1/4 & \quad w = \pm \beta\,, \\
(w^2 -\beta^2)^{-1/2}\,J_1\big(\sqrt{w^2 - \beta^2}\big) &\quad w \in \mathbb R \setminus [-\beta,\,\beta]\,.
\end{array} \right.
\end{equation}
By the power series expansion of the modified Bessel function $I_1$, we obtain for $w \in (- \beta,\,\beta)$,
$$
(\beta^2 - w^2)^{-1/2}\,I_1\big(\sqrt{\beta^2 - w^2}\big) = \frac{1}{2}\,\sum_{k=0}^{\infty} \frac{1}{4^k\, k!\, (k+1)!} \,(\beta^2 - w^2)^k
$$
such that ${\hat \varphi}_{\sinh}(v)$ is decreasing for $v \in \big[0,\,N_1\big(1 - \frac{1}{2\sigma}\big)\big]$. The frequency $v = \frac{N}{2}$ corresponds to the scaled frequency  $w = \frac{\pi m}{\sigma} < \beta = 2\pi m \big(1 - \frac{1}{2\sigma}\big)$.
Hence we have
\begin{eqnarray*}
\min_{n \in I_N} {\hat \varphi}_{\sinh}(n) &=& {\hat \varphi}_{\sinh}\big(\frac{N}{2}\big)\\
&=& \frac{m \pi \beta}{N_1 \sinh \beta}\,\big(\beta^2 - \frac{\pi^2 m^2}{\sigma^2}\big)^{-1/2} \,I_1\big(\sqrt{\beta^2 - \frac{\pi^2 m^2}{\sigma^2}}\big)\\
&=& \frac{\beta}{2 N_1 \sinh \beta}\,\big(1 - \frac{1}{\sigma}\big)^{-1/2} \,I_1\big(2 \pi m \sqrt{1 - \frac{1}{\sigma}}\big)\,.
\end{eqnarray*}
By $m\ge 2$ and $\sigma \ge \frac{5}{4}$, we maintain
$$
2 \pi m \sqrt{1 - \frac{1}{\sigma}} \ge 4 \pi\,\sqrt{1 - \frac{1}{\sigma}} \ge x_0 := \frac{4 \pi}{\sqrt 5}\,.
$$
By the inequality (see \cite{Bar} or \cite[Lemma 3.3]{PoTa20})
$$
\sqrt{2 \pi x_0}\, {\mathrm e}^{-x_0}\,I_1(x_0) \le \sqrt{2 \pi x}\, {\mathrm e}^{-x}\,I_1(x)\,, \quad x \ge x_0\,,
$$
we preserve for all $x \ge x_0$,
$$
I_1(x) \ge \sqrt x_0\, {\mathrm e}^{-x_0}\,I_1(x_0)\, x^{-1/2}\, {\mathrm e}^x > \frac{2}{5}\,x^{-1/2}\, {\mathrm e}^x
$$
and hence
\begin{equation}
\label{eq:hatvarphisinhN/2}
{\hat \varphi}_{\sinh}\big(\frac{N}{2}\big) \ge \frac{\beta}{5 \,N_1\,\sqrt{2 \pi m}\, \sinh \beta}\,\big(1 - \frac{1}{\sigma}\big)^{-3/4}\,{\mathrm e}^{2 \pi m \sqrt{1 - 1/\sigma}}\,.
\end{equation}
Now we estimate ${\hat \varphi}_{\sinh}\big(\frac{N_1 w}{2\pi m}\big)$ for $|w| > \beta =b m$.
For $|w| > \beta$, $\sigma \ge \frac{5}{4}$, and $N\ge 8$, it holds
$$
\frac{N_1 |w|}{2 \pi m} > N_1 \big(1 - \frac{1}{2 \sigma}\big) = N\, \big(\sigma - \frac{1}{2}\big) \ge \frac{3 N}{4} \ge 6\,.
$$
By an inequality of the Bessel function $J_1$ (see \cite{Kra} or \cite[Lemma 3.2]{PoTa20}), we have for all $x \ge 6$,
$$
|J_1(x)| \le \frac{1}{\sqrt x}\,.
$$
Using \eqref{eq:FTsinhwindow}, it follows that for $|w| > \beta$,
\begin{eqnarray*}
\big|{\hat \varphi}_{\sinh}\big(\frac{N_1 w}{2\pi m}\big)\big| &=& \frac{\beta  \pi m}{N_1\,\sinh \beta}\,(w^2 - \beta^2)^{-1/2}\,\big|J_1\big(\sqrt{w^2 - \beta^2}\big)\big|\\
&\le& \frac{\beta  \pi m}{N_1\,\sinh \beta}\,(w^2 - \beta^2)^{-3/4}
\end{eqnarray*}
i.e., for $|v| > N_1 \big(1 - \frac{1}{2\sigma}\big)$,
\begin{equation}
\label{eq:hatvarphisinh}
|{\hat \varphi}_{\sinh}(v)| \le \frac{ \sqrt{\pi m} }{{\sqrt 2}\,N_1\, \sinh \beta}\,\Big(\frac{v^2}{N_1^2} - \big(1 - \frac{1}{2 \sigma}\big)^2\Big)^{-3/4} \big(1 - \frac{1}{2 \sigma}\big)\,.
\end{equation}
For each $n \in I_N$ with $|n \pm N_1| > N_1 \big(1 - \frac{1}{2\sigma}\big)$, we obtain by \eqref{eq:FTsinhwindow} and $|J_1(x)| \le \frac{1}{2}\,|x|$ that
\begin{equation}
\label{eq:FTsinh(n+N1)}
|{\hat \varphi}_{\sinh}(n \pm N_1)| \le \frac{\pi m \beta }{2\,N_1\, \sinh \beta}\,.
\end{equation}
In the case $- \frac{N}{2} + N_1 = N_1 \big(1 - \frac{1}{2 \sigma}\big)$, we get by \eqref{eq:FTsinhwindow} that
$$
{\hat \varphi}_{\sinh}\big(- \frac{N}{2} + N_1 \big) = \frac{\pi m \beta}{4\,N_1 \sinh \beta}\,,
$$
i.e., we can use \eqref{eq:FTsinh(n+N1)} for all $n \in I_N$.

For all $n \in I_N$ and $r \in {\mathbb Z}\setminus \{0,\, \pm 1\}$, we have by $\sigma \ge \frac{5}{4}$ that
$$
\big|\frac{n}{N_1} + r| \ge 2 - \frac{1}{2\sigma} \ge 2 - \frac{2}{5} = \frac{8}{5}\,.
$$
Since the decreasing function $h: \big[\frac{8}{5},\,\infty\big) \to \mathbb R$,
$$
h(x):=x^{3/2}\,\Big(x^2 - \big(1 - \frac{1}{2 \sigma}\big)^2\Big)^{-3/4} = \Big(1 - \big(1 - \frac{1}{2 \sigma}\big)^2\, x^{-2}\Big)^{-3/4}\,,
$$
is bounded by
$$
h\big(\frac{8}{5}\big) = \Big(1 - \frac{25}{64}\,\big(1 - \frac{1}{2\sigma}\big)^2\Big)^{-3/4} < \big(1 - \frac{25}{64}\big)^{-3/4} < \frac{3}{2}\,,
$$
we receive for all  $n\in I_N$ and $r \in \mathbb Z \setminus \{0,\,\pm 1\}$,
\begin{eqnarray*}
\big|{\hat \varphi}_{\sinh}(n + r N_1)\big| &\le& \frac{\sqrt{\pi m}}{{\sqrt 2}\,N_1\, \sinh \beta} \Big(\big(r + \frac{n}{N_1}\big)^2 - \big(1 - \frac{1}{2 \sigma}\big)^2\Big)^{-3/4}\\
&\le& \frac{3\, \sqrt {\pi m}}{2{\sqrt 2}\, N_1\, \sinh \beta} \big|r + \frac{n}{ N_1} \big|^{-3/2} \,.
\end{eqnarray*}
Therefore we get by \eqref{eq:sumxmu} that for all $n \in I_N$,
\begin{eqnarray*}
& &\sum_{r\in \mathbb Z \setminus \{0,\,\pm 1\}} \big|{\hat \varphi}_{\mathrm{Bessel}}(n + r N_1)\big| \le \frac{3\, \sqrt {\pi m}}{2{\sqrt 2}\, N_1\, \sinh \beta} \sum_{r\in \mathbb Z \setminus \{0,\,\pm 1\}} \big|r + \frac{n}{ N_1} \big|^{-3/2}\\
& & \le \frac{3\, \sqrt {2 \pi m}}{N_1\, \sinh \beta}\,\big(1 - \frac{1}{2 \sigma}\big)^{-1/2}
\end{eqnarray*}
and hence by \eqref{eq:FTsinh(n+N1)},
\begin{equation}
\label{eq:sumFTsinh}
\sum_{r\in \mathbb Z \setminus \{0\}} \big|{\hat \varphi}_{\mathrm{Bessel}}(n + r N_1)\big| \le \frac{1}{N_1 \sinh \beta}\Big[\pi m \beta + 3 \sqrt{2\pi m}\, \big(1 - \frac{1}{2 \sigma}\big)^{-1/2}\Big]\,.
\end{equation}
Then from \eqref{eq:FTsinhwindow} and \eqref{eq:FTsinh(n+N1)} it follows that
\begin{eqnarray*}
e_{\sigma,N}(\varphi_{\sinh}) &\le& \frac{1}{{\hat \varphi}_{\sinh}(N/2)} \max_{n \in I_N} \sum_{r\in \mathbb Z \setminus \{0\}} \big|{\hat \varphi}_{\mathrm{Bessel}}(n + r N_1)\big|\\
&\le& \Big[5 \pi m \sqrt{2 \pi m} + 3 \big(1 - \frac{1}{2 \sigma}\big)^{-3/2}\Big]\,\big(1 - \frac{1}{\sigma}\big)^{3/4}\,{\mathrm e}^{-2 \pi m\, \sqrt{1 - 1/\sigma}}\,.
\end{eqnarray*}
Using $5\pi \sqrt{2 \pi} < 40$, we summarize:

\begin{Theorem}
\label{Thm:esigmaNsinh}
Let $N \in 2 \mathbb N$, $N\ge 8$, and $\sigma \ge \frac{5}{4}$ be given, where $N_1 = \sigma N \in 2 \mathbb N$. Further let $m \in \mathbb N \setminus \{1\}$ with $2m \ll N_1$ and $\beta = 2 \pi m \big(1 - \frac{1}{2\sigma}\big)$.

Then the $C(\mathbb T)$-error constant of the $\sinh$-type window function \eqref{eq:sinhwindow} can be estimated by
$$
e_{\sigma}({\varphi}_{\sinh}) \le \Big[ 40\, m^{3/2} + 3 \big(1 - \frac{1}{2\sigma}\big)^{-3/2}\Big]\, \big(1 - \frac{1}{\sigma}\big)^{3/4} \,{\mathrm e}^{-2 \pi m\, \sqrt{1 - 1/\sigma}}\,,
$$
i.e., the $\sinh$-type window function \eqref{eq:sinhwindow} is convenient for ${\mathrm{NFFT}}$.
\end{Theorem}

For a fixed oversampling factor $\sigma \ge \frac{5}{4}$, the $C(\mathbb T)$-error constant of \eqref{eq:sinhwindow} decays exponentially with the truncation parameter $m \ge 2$. On the other hand, the
computational cost of NFFT increases with respect to $m$ (see \cite[pp.~380--381]{PlPoStTa18}) such that $m$ should be not too large. For $\sigma = 2$ and $m = 4$, we obtain $e_{\sigma}({\varphi}_{\sinh}) \le 3.7\cdot 10^{-6}$.

\section{Continuous exp-type window function}\label{Sec:expWindow}

For fixed shape parameter $\beta = b m$ with $m\in \mathbb N \setminus \{1\}$, $b = 2 \pi\,\big(1 - \frac{1}{2 \sigma}\big)$, and oversampling factor $\sigma \ge \frac{5}{4}$, we consider the \emph{continuous} $\exp$-\emph{type window function}
\begin{equation}
\label{eq:expwindow}
\varphi_{\mathrm{cexp}}(x) := \left\{ \begin{array}{ll}
\frac{1}{{\mathrm e}^{\beta} - 1} \Big({\mathrm e}^{\beta \,\sqrt{1 - (N_1x/m)^2}} - 1\Big) & \quad x \in I\,, \\ [1ex]
0 & \quad x \in \mathbb R \setminus  I\,.
\end{array}\right.
\end{equation}
Obviously, $\varphi_{\mathrm{cexp}} \in  \Phi_{m,N_1}$ is a continuous window function. Note that a discontinuous version of this window function was suggested in \cite{BaMaKl18, Ba20}. A corresponding error estimate for the NFFT was proved in \cite{Ba20}, where an asymptotic value of its Fourier transform was determined for $\beta \to \infty$ by saddle point integration. We present new explicit error estimates for fixed shape parameter $\beta$ of moderate size.
\medskip

In the following, we present a new approach to an error estimate for the NFFT with
the continuous $\exp$-type window function \eqref{eq:expwindow}.
Unfortunately, the Fourier transform of \eqref{eq:expwindow} is unknown analytically. Therefore we represent \eqref{eq:expwindow} as sum
\begin{equation}
\label{eq:splitvarphicexp}
\varphi_{\mathrm{cexp}}(x) = \psi(x) + \rho(x)\,,
\end{equation}
where the Fourier transform of $\psi$ is known and where the correction term $\rho$ has small magnitude $|\rho|$. We choose
$$
\psi(x) := \left\{ \begin{array}{ll}
\frac{2}{{\mathrm e}^{\beta} - 1} \, \sinh \big(\beta \,\sqrt{1 - (N_1 x)^2/m^2}\big) & \quad x \in I\,, \\
0 & \quad x \in \mathbb R \setminus  I
\end{array}\right.
$$
and
$$
\rho(x) := \left\{ \begin{array}{ll}
\frac{1}{{\mathrm e}^{\beta} - 1} \, \Big({\mathrm e}^{-\beta \,\sqrt{1 - (N_1 x)^2/m^2}}-1\Big) & \quad x \in I\,, \\ [1ex]
0 & \quad x \in \mathbb R \setminus  I\,.
\end{array}\right.
$$
The Fourier transform of $\psi$ reads as follows (see \cite[p.~38]{Oberh90})
\begin{equation}
\label{eq:FTpsi}
{\hat \psi}(v) = \frac{2 \pi m \beta}{({\mathrm e}^{\beta} - 1)\,N_1} \left\{ \begin{array}{ll}
\frac{1}{2 \pi m}\big(\big(1 - \frac{1}{2 \sigma}\big)^2 - \frac{v^2}{N_1^2}\big)^{-1/2} \,I_1 \big( 2\pi m \sqrt{\big(1 - \frac{1}{2 \sigma}\big)^2 - \frac{v^2}{N_1^2}}\big) & \quad |v| < N_1 \big(1 - \frac{1}{2\sigma}\big)\,,\\
\frac{1}{4}  & \quad v = \pm N_1 \big(1 - \frac{1}{2\sigma}\big)\,,\\
\frac{1}{2 \pi m}\big( \frac{v^2}{N_1^2} - \big(1 - \frac{1}{2 \sigma}\big)^2\big)^{-1/2} \,J_1 \big( 2\pi m \sqrt{ \frac{v^2}{N_1^2} - \big(1 - \frac{1}{2 \sigma}\big)^2}\big) & \quad |v| > N_1 \big(1 - \frac{1}{2\sigma}\big)\,.
\end{array} \right.
\end{equation}
Since $\rho$ is even and $\rho  |_{[0,\, \frac{m}{N_1}]}$ is increasing, we have
$
0 \ge \rho(x) \ge \rho(0) = \frac{1}{{\mathrm e}^{\beta} - 1}\,({\mathrm e}^{-\beta}-1) = -{\mathrm e}^{-\beta}\,.
$

\begin{table}[t]
\begin{center}
\begin{tabular}{c|c}
$m$ & $-\rho(0)= {\mathrm e}^{-\beta}$ \\ \hline
$2$  & $8.06 \cdot 10^{-5}$\\
$3$  & $7.24 \cdot 10^{-7}$\\
$4$  & $6.51 \cdot 10^{-9}$\\
$5$  & $5.85 \cdot 10^{-11}$\\
$6$  & $5.25 \cdot 10^{-13}$
\end{tabular}
\end{center}
\caption{Upper bounds of the function $-\rho$ for $m=2,\,3,\,\ldots,\,6$ and $\sigma = 2$.}
\end{table}

Since $\rho$ has small absolute values in the small support $\bar I$, the Fourier transform $\hat \rho$ is small too and it holds
$$
|{\hat \rho}(v)| = \big|\int_I \rho(x)\,{\mathrm e}^{- 2\pi {\mathrm i}\,v x}\,{\mathrm d}x \big| \le \frac{2 m}{N_1}\,{\mathrm e}^{-\beta}\,.
$$
Substituting $t= N_1 x/m$, we determine the Fourier transform
\begin{eqnarray*}
{\hat \varphi}_{\mathrm{cexp}}(v) &=& \int_I \varphi_{\mathrm{cexp}}(x)\,{\mathrm e}^{-2\pi {\mathrm i}\,v x}\, {\mathrm d}x
= \int_I \psi(x)\,{\mathrm e}^{-2\pi {\mathrm i}\,v x}\, {\mathrm d}x + \int_I \rho(x)\,{\mathrm e}^{-2\pi {\mathrm i}\,v x}\, {\mathrm d}x \\
&=& \frac{2 m}{({\mathrm e}^{\beta} - 1)\,N_1}\Big[2 \,\int_0^1 \sinh\big(\beta \sqrt{1-t^2}\big)\,\cos\frac{2\pi m v t}{N_1}\, {\mathrm d}t\\
& &+\, \int_0^1 \big({\mathrm e}^{-\beta \sqrt{1 - t^2}}-1\big) \, \cos\frac{2\pi m v t}{N_1}\,{\mathrm d}t\Big]\,.
\end{eqnarray*}
For simplicity, we introduce the scaled frequency  $w:= 2 \pi m v/N_1$ such that
\begin{eqnarray}
\label{eq:FTexptypewindow}
{\hat \varphi}_{\mathrm{cexp}}\big(\frac{N_1 w}{2 \pi m}\big) &=& \frac{2 m}{({\mathrm e}^{\beta}- 1)\,N_1} \Big[2\,\int_0^1 \sinh\big(\beta \sqrt{1-t^2}\big)\,\cos(w t)\, {\mathrm d}t\nonumber \\
& &+\, \int_0^1 \big({\mathrm e}^{-\beta \sqrt{1 - t^2}}-1\big) \, \cos(wt)\,{\mathrm d}t\Big]\,.
\end{eqnarray}
From \cite[p. 38]{Oberh90} it follows that
$$
\int_0^1 \sinh\big(\beta \sqrt{1-t^2}\big)\,\cos(w t)\, {\mathrm d}t = \frac{\pi \beta}{2} \left\{\begin{array}{ll}
(\beta^2 - w^2)^{-1/2}\,I_1\big(\sqrt{\beta^2 - w^2}\big) & \quad w \in (-\beta,\, \beta)\,,\\
1/4 & \quad w = \pm \beta\,, \\
(w^2 - \beta^2)^{-1/2}\,J_1\big(\sqrt{w^2 - \beta^2}\big) & \quad w \in \mathbb R \setminus [-\beta,\, \beta]\,,
\end{array} \right.
$$
where $I_1$ denotes the modified Bessel function and $J_1$ the Bessel function of first order. Therefore we consider
\begin{eqnarray}
\hat \rho\big(\frac{N_1 w}{2 \pi m}\big)&=&   \frac{2 m}{({\mathrm e}^{\beta}- 1)\,N_1} \,\int_0^1 \big({\mathrm e}^{-\beta \sqrt{1 - t^2}}-1\big) \, \cos(wt)\,{\mathrm d}t \nonumber \\
&=&  \frac{m}{({\mathrm e}^{\beta}- 1)\,N_1} \,\int_{-1}^1 \big({\mathrm e}^{-\beta \sqrt{1 - t^2}}-1\big) \, {\mathrm e}^{{\mathrm i}\,wt}\,{\mathrm d}t \label{eq:hatrhoinw}
\end{eqnarray}
as correction term of \eqref{eq:FTexptypewindow}.
Now we estimate the integral
\begin{equation}
\label{eq:mathcalI(w)}
{\mathcal I}(w):= \int_{-1}^1 \big({\mathrm e}^{-\beta \sqrt{1 - t^2}}-1\big) \, {\mathrm e}^{{\mathrm i}\,wt}\,{\mathrm d}t
\end{equation}
by complex contour integrals.

\begin{Lemma}
\label{Lemma:est|I(w)|}
For each $w\in \mathbb R$ with $|w| \ge \beta = b m$, $m \in \mathbb N \setminus \{1\}$, we have
$$
|{\mathcal I}(w)|\le  \big(2 + \frac{2}{\mathrm e}\big)\,\beta\,\sqrt[4]5\,|w|^{-5/4} + 4\,|w|^{-1}\,{\mathrm e}^{-\sqrt{|w|}} + \big(2\,{\mathrm e}^{-\sqrt 2 \beta} + 1\big)\,{\mathrm e}^{-|w|}\,.
$$
\end{Lemma}

\emph{Proof}. Here we use the same technique as in \cite[Lemma 10]{Ba20}, where the integral
$$
\int_{-1}^1 \big({\mathrm e}^{\beta \sqrt{1 - t^2}}-1\big) \, {\mathrm e}^{{\mathrm i}\,wt}\,{\mathrm d}t
$$
for $|w| \ge \beta^4$ was estimated.
Since ${\mathcal I}(w) = {\mathcal I}(-w)$, we consider only the case $w \ge \beta$.
Let $C_1$ be the line segment from $-1$ to $-1 + {\mathrm i}$. Further, $C_2$ is the line segment from $-1 + {\mathrm i}$ to $1 + {\mathrm i}$, $C_3$ is the line
segment from $1 + {\mathrm i}$ to $1$, and $C_4$ is the line segment from 1 to $-1$.
Since the principal square root function is holomorphic in $\mathbb C$ except the nonpositive real axis, the complex function
$$
g(z) :=\big({\mathrm e}^{-\beta \sqrt{1 - z^2}}-1\big) \, {\mathrm e}^{{\mathrm i}\,w z}
$$
is holomorphic in $\mathbb C$ except the set $(- \infty,\,-1] \cup [1,\,\infty)$. Hence $g$ is holomorphic on the interior of the closed curve $C := C_1 \cup C_2 \cup C_3 \cup C_4$ and continuous in $D \cup C$,
since simple calculation shows that
$$
\lim_{\stackrel{z \to -1}{z\in D}} g(z) = \lim_{\stackrel{z \to 1}{z \in D}} g(z) = 0\,.
$$
Then the stronger form of Cauchy's Integral Theorem (see \cite{Be43}) provides
\begin{equation}
\label{eq:Cauchy}
{\mathcal I}(w) = {\mathcal I}_1(w) + {\mathcal I}_2(w) + {\mathcal I}_3(w)
\end{equation}
with the contour integrals
$$
{\mathcal I}_k(w) =\int_{C_k} \big({\mathrm e}^{-\beta \sqrt{1 - z^2}}-1\big) \, {\mathrm e}^{{\mathrm i}\,w z}\,{\mathrm d}z\,, \quad k = 1,\,2,\,3\,.
$$
Note that ${\mathcal I}_3(w)$ is the complex conjugate of ${\mathcal I}_1(w)$ such that $|{\mathcal I}_3(w)| = |{\mathcal I}_1(w)|$.

The line segment $C_2$ can be parametrized by $z = t + {\mathrm i}$, $t\in[-1,\,1]$ such that
$$
{\mathcal I}_2(w) = {\mathrm e}^{-w}\,\int_{-1}^1 \big({\mathrm e}^{-\beta \sqrt{2 - t^2 - 2{\mathrm i}t}} - 1\big)\, {\mathrm e}^{{\mathrm i}\,w t}\,{\mathrm d}t
$$
and hence
$$
|{\mathcal I}_2(w)| \le  {\mathrm e}^{-w}\,\int_{-1}^1 \big|{\mathrm e}^{-\beta \sqrt{2 - t^2 - 2{\mathrm i}t}}\big| \,{\mathrm d}t + 2\,{\mathrm e}^{-w} \,.
$$
Then we have
$
\big|{\mathrm e}^{-\beta \sqrt{2 - t^2 - 2{\mathrm i}t}}\big| = {\mathrm e}^{-\beta\, {\mathrm{Re}}\sqrt{2 - t^2 - 2{\mathrm i}t}}\,.
$
Since for $t\in [-1,\,1]$,
$$
|2 - t^2 -2 {\mathrm i}\, t| = \sqrt{t^2 + 4} \in [2,\, \sqrt 5]\,,
$$
we obtain the estimate
$\sqrt 2 \le {\mathrm{Re}}\sqrt{2 - t^2 - 2{\mathrm i}t} \le \sqrt[4]5$ for $t\in [-1,\,1]$.

Note that ${\mathrm{Re}}\sqrt{2 - t^2 - 2{\mathrm i}t}> 0$ for all $t\in [-1,\,1]$. Thus we have
$$
|{\mathcal I}_2(w)| \le 2 \,{\mathrm e}^{-w}\,\big({\mathrm e}^{- {\sqrt 2} \beta} + 1\big)\,.
$$
A parametrization of the line segment $C_1$ is $z = -1 + {\mathrm i}\,t$, $t\in [0,\,1]$, such that
$$
{\mathcal I}_1(w) = {\mathrm i}\,{\mathrm e}^{- {\mathrm i}w}\, \int_0^1 \big({\mathrm e}^{-\beta\,\sqrt{2 {\mathrm i}t + t^2}}-1\big)\,{\mathrm e}^{-w t}\,{\mathrm d}t\,.
$$
For $w \ge \beta >0$, we split ${\mathcal I}_1(w)$ into the sum of two integrals
\begin{eqnarray*}
{\mathcal I}_{1,0}(w) &:=& {\mathrm i}\,{\mathrm e}^{- {\mathrm i}w}\,\int_0^{w^{-1/2}}\big({\mathrm e}^{-\beta\,\sqrt{2 {\mathrm i}t + t^2}}-1\big)\,{\mathrm e}^{-w t}\,{\mathrm d}t\,,\\ [1ex]
{\mathcal I}_{1,1}(w) &:=& {\mathrm i}\,{\mathrm e}^{- {\mathrm i}w}\,\int_{w^{-1/2}}^1\big({\mathrm e}^{-\beta\,\sqrt{2 {\mathrm i}t + t^2}}-1\big)\,{\mathrm e}^{-w t}\,{\mathrm d}t\,.
\end{eqnarray*}
Since
$
2 \le |\sqrt{2 {\mathrm i} + t}| \le \sqrt[4]5$, $t \in [0,\,1]\,,
$
the integral ${\mathcal I}_{1,0}(w)$ is bounded in magnitude by
$$
|{\mathcal I}_{1,0}(w)| \le \max_{t\in [0,\,w^{-1/2}]} \big(1 - {\mathrm e}^{-\beta \sqrt{t}\,\sqrt[4]5}\big)\,\int_0^{w^{-1/2}}{\mathrm e}^{-w t}\,{\mathrm d}t\,.
$$
From
$$
\int_0^{w^{-1/2}}{\mathrm e}^{-w t}\,{\mathrm d}t = \frac{1 - {\mathrm e}^{-\sqrt w}}{w}\le \frac{1}{w}
$$
and
$$
\max_{t\in [0,\,w^{-1/2}]} \big(1 - {\mathrm e}^{-\beta \sqrt{t}\,\sqrt[4]5}\big) = 1 - {\mathrm e}^{-\beta\,\sqrt[4]{5/w}} \le \beta\,\sqrt[4]{\frac{5}{w}}
$$
it follows that
$$
|{\mathcal I}_{1,0}(w)| \le \big(1 + {\mathrm e}^{-1}\big)\,\beta \sqrt[4]{5}\,w^{-5/4}\,.
$$
Above we have used the simple inequality $1 - {\mathrm e}^{-x} \le x$ for $x \ge 0$.

Finally we estimate the integral ${\mathcal I}_{1,1}(w)$ as follows
$$
|{\mathcal I}_{1,1}(w)| \le \int_{w^{-1/2}}^1 {\mathrm e}^{- \beta {\sqrt t}\,{\mathrm{Re}}\sqrt{2\mathrm i + t} - wt}\, {\mathrm d}t + \int_{w^{-1/2}}^1 {\mathrm e}^{ - w t}\,{\mathrm d}t\,.
$$
From ${\mathrm{Re}}\sqrt{2\mathrm i + t}> 0$ for all $t\in [0,\,1]$ and
$$
\int_{w^{-1/2}}^1 {\mathrm e}^{ - w t}\,{\mathrm d}t = w^{-1}\,\big({\mathrm e}^{-\sqrt w} - {\mathrm e}^{-w}\big)\le w^{-1}\,{\mathrm e}^{-\sqrt w}
$$
it follows that
$$
|{\mathcal I}_{1,1}(w)| \le \int_{w^{-1/2}}^1 {\mathrm e}^{ - w t}\,{\mathrm d}t +  w^{-1}\,{\mathrm e}^{-\sqrt w} \le 2\, w^{-1}\,{\mathrm e}^{-\sqrt w}\,.
$$
Thus we receive for $w\ge \beta$,
\begin{eqnarray*}
|{\mathcal I}(w)| &\le& |{\mathcal I}_1(w)| + |{\mathcal I}_2(w)| + |{\mathcal I}_3(w)| = 2\,|{\mathcal I}_1(w)| + |{\mathcal I}_2(w)|\\
&\le& 2\,|{\mathcal I}_{1,0}(w)| + 2\,|{\mathcal I}_{1,1}(w)| + |{\mathcal I}_2(w)|\\
&\le&  \big(2 + \frac{2}{\mathrm e}\big)\,\beta\,\sqrt[4]5\,w^{-5/4} + 4\,w^{-1}\,{\mathrm e}^{-\sqrt w} + \big(2\,{\mathrm e}^{-{\sqrt 2} \beta} + 1\big)\,{\mathrm e}^{-w}\,.
\end{eqnarray*}
This completes the proof. \qedsymbol

\begin{Remark} In the proof of  Lemma \ref{Lemma:est|I(w)|}, it was shown that for real $w$ with $|w| \ge \beta$, the contour integral ${\mathcal I}_1(w)$ can be
estimated by
$$
|{\mathcal I}_1(w)| \le \big(1 + \frac{1}{\mathrm e}\big)\,\beta \sqrt[4]5\,|w|^{-5/4} + 2\,|w|^{-1}\,{\mathrm e}^{-\sqrt{|w|}}\,.
$$
Splitting the integral ${\mathcal I}_1(w)$ in the form
$$
{\mathcal I}_1(w) ={\mathrm i}\,{\mathrm e}^{- {\mathrm i}\,w}\,\Big(\int_0^{(\log \sqrt w)/w} + \int_{(\log \sqrt w)/w}^1\Big)\big({\mathrm e}^{- \beta\,\sqrt{2 {\mathrm i}\,t + t^2}} - 1\big)\,{\mathrm e}^{-w t}\,{\mathrm d}t\,,
$$
this yields the better estimate for $|w| \ge \beta$,
$$
|{\mathcal I}_1(w)| \le \big(\frac{1}{\sqrt 2}\,\beta \sqrt[4]5 \,\sqrt{\log |w|} + 2\big)\,|w|^{-3/2}\,.
$$
\end{Remark}

\begin{Lemma}
\label{Lemma:estint01bysinc} Let $\beta = bm$ be given.
Then for each $w \in \mathbb R$, it holds the estimate
$$
- {\mathrm{sinc}}\,w - \gamma(m,\sigma) \le \int_0^1 \big({\mathrm e}^{-\beta \sqrt{1 - t^2}} - 1\big) \, \cos(wt)\,{\mathrm d}t \le - {\mathrm{sinc}}\,w + \gamma(m,\sigma)
$$
with the small positive constant
$$
\gamma(m,\sigma) := \int_0^1 {\mathrm e}^{-\beta \sqrt{1 - t^2}}\,{\mathrm d}t\,.
$$
\end{Lemma}

\emph{Proof}.
For each $w \in \mathbb R$, we have
\begin{eqnarray*}
& &\big|\int_0^1 \big({\mathrm e}^{-\beta \sqrt{1 - t^2}}-1\big) \, \cos(wt)\,{\mathrm d}t + \mathrm{sinc}\,w\big| = \big|\int_0^1 {\mathrm e}^{-\beta \sqrt{1 - t^2}} \, \cos(wt)\,{\mathrm d}t\big|\\
& & \le \, \int_0^1 {\mathrm e}^{-\beta \sqrt{1 - t^2}}\,{\mathrm d}t = \gamma(m,\sigma)
\end{eqnarray*}
such that
$$
\max_{w \in \mathbb R} \big| \int_0^1 \big({\mathrm e}^{-\beta \sqrt{1 - t^2}}-1\big) \, \cos(wt)\,{\mathrm d}t + \mathrm{sinc}\,w\big| \le \gamma(m, \sigma)\, ,
$$
see Table 6.2. \qedsymbol

\begin{table}[t]
\begin{center}
\begin{tabular}{c|c}
$m$ & $\gamma(m, 2)$ \\ \hline
$2$  & $1.17 \cdot 10^{-2}$\\
$3$  & $5.08 \cdot 10^{-3}$\\
$4$  & $2.84 \cdot 10^{-3}$\\
$5$  & $1.81 \cdot 10^{-3}$\\
$6$  & $1.25 \cdot 10^{-3}$
\end{tabular}
\end{center}
\caption{Maximum error $\gamma(m, \sigma)$ between $\int_0^1 \big({\mathrm e}^{-\beta \sqrt{1 - t^2}}-1\big) \, \cos(wt)\,{\mathrm d}t$ and $-\mathrm{sinc}\,w$ for $m=2,\,3,\,\ldots,\,6$ and $\sigma$ = 2.}
\end{table}
\medskip

Using \eqref{eq:FTexptypewindow} and Lemma \ref{Lemma:estint01bysinc}, we receive for $w \in (-\beta,\,\beta)$,
\begin{eqnarray}
\label{eq:FTexpwindowinw}
{\hat \varphi}_{\mathrm{cexp}}\big(\frac{N_1 w}{2 \pi m}\big)&\ge& \frac{2 m}{({\mathrm e}^{\beta} - 1)\,N_1}\Big[\pi \beta\,(\beta^2 - w^2)^{-1/2}\,I_1\big(\sqrt{\beta^2 - w^2}\big)\\
& & -\, \mathrm{sinc}\,w - \gamma(m,\sigma)\Big]\,.
\end{eqnarray}
The function $h:\,[0, \beta] \to \mathbb R$ defined by
\begin{eqnarray*}
h(w) &:=& \pi \beta\,(\beta^2 - w^2)^{-1/2}\,I_1\big(\sqrt{\beta^2 - w^2}\big) - \mathrm{sinc}\,w \\
&=& \frac{\pi \beta}{2} \sum_{s = 0}^{\infty} \frac{(\beta^2 - w^2)^s}{2^{2s}\,s!\,(s+1)!}  - \mathrm{sinc}\,w
\end{eqnarray*}
has the derivative
\begin{equation}
\label{eq:h'(w)}
h'(w) = - \pi \beta\,\sum_{s = 1}^{\infty}  \frac{(\beta^2 - w^2)^{s-1}}{2^{2s}\,(s-1)!\,(s+1)!}  - \frac{{\mathrm d}}{{\mathrm d}w}\,\mathrm{sinc}\,w\,.
\end{equation}
Obviously, it holds $h'(0) = 0$, and $h'(w) < 0$ for $w \in (0,\,\pi]$. For $w \in [\pi,\, \beta]$ and $m \ge 2$, we obtain
$$
\big|\frac{{\mathrm d}}{{\mathrm d}w}\,\mathrm{sinc}\,w\big| \le \frac{1}{\pi^2} + \frac{1}{\pi}
$$
and hence by \eqref{eq:h'(w)},
$$
h'(w) < - \frac{\pi \beta}{8} + \frac{1}{\pi^2} + \frac{1}{\pi} < 0\,.
$$
Thus  $h(w)$ is decreasing with respect to  $w \in [0,\, \beta]$.

For $\beta = bm$ and $v = \frac{N_1 w}{2\pi m}\in \big(-N_1\big(1- \frac{1}{2\sigma}\big),\,N_1 \big(1- \frac{1}{2\sigma}\big)\big)$, the inequality \eqref{eq:FTexpwindowinw} implies
\begin{eqnarray*}
{\hat \varphi}_{\mathrm{cexp}}(v) &\ge& \frac{2 m}{({\mathrm e}^{\beta} - 1)\,N_1}\Big[2\, \Big(\big(1 - \frac{1}{2\sigma}\big)^2 - \frac{v^2}{N_1^2}\Big)^{-1/2}\,I_1\Big(2\pi m \sqrt{\big(1 - \frac{1}{2\sigma}\big)^2 - \frac{v^2}{N_1^2}}\Big)\\
& & -\, \mathrm{sinc}\frac{2 \pi m v}{N_1}- \gamma(m,\sigma)\Big]\,.
\end{eqnarray*}
Consequently,  the Fourier transform ${\hat \varphi}_{\mathrm{cexp}}(v)$ is positive and decreasing for
$$
v \in \big[0,\, N_1 \big(1- \frac{1}{2\sigma}\big)\big)\,.
$$

Hence we obtain
\begin{eqnarray*}
\min_{n\in I_N} {\hat \varphi}_{\mathrm{cexp}}(n) &=& {\hat \varphi}_{\mathrm{cexp}}\big(\frac{N}{2}\big)\\
&\ge& \frac{2 m}{({\mathrm e}^{\beta} - 1)\,N_1}\Big[\frac{b}{2}\,\big(1 - \frac{1}{\sigma}\big)^{-1/2}\,I_1\big(2 \pi m \sqrt{1 - \frac{1}{\sigma}}\big) - \mathrm{sinc}\frac{\pi m}{\sigma} - \gamma(m,\sigma)\Big]\\
&\ge& \frac{2 m}{({\mathrm e}^{\beta} - 1)\,N_1}\Big[\frac{b}{2}\,\big(1 - \frac{1}{\sigma}\big)^{-1/2}\,I_1\big(2 \pi m \sqrt{1 - \frac{1}{\sigma}}\big) - 1 - \gamma(m,\sigma)\Big]\,.
\end{eqnarray*}
From $m \ge 2$ and $\sigma \ge \frac{5}{4}$  it follows that
$$
2 \pi m \sqrt{1 - \frac{1}{\sigma}} \ge 4 \pi \,\sqrt{1 - \frac{1}{\sigma}}\ge x_0 := \frac{4 \pi}{\sqrt 5}\,.
$$

Hence by the inequality (see \cite{Bar} or \cite[Lemma 3.3]{PoTa20})
$$
\sqrt{2 \pi x_0}\,{\mathrm e}^{-x_0}\, I_1(x_0) \le \sqrt{2 \pi x}\,{\mathrm e}^{-x}\, I_1(x)\,, \quad x \ge x_0\,,
$$
we sustain that
$$
I_1(x) \ge \sqrt x_0 \,{\mathrm e}^{-x_0}\,I_1(x_0)\,x^{-1/2}\,{\mathrm e}^x\,, \quad x \ge x_0\,.
$$
Thus for $x = 2 \pi m \sqrt{1- 1/\sigma}$, we get the estimate
$$
I_1\big( 2 \pi m \sqrt{1 - \frac{1}{\sigma}}\big) \ge \frac{2}{5\,\sqrt{2 \pi m}} \, \big(1 - \frac{1}{\sigma}\big)^{-1/4}\,{\mathrm e}^{2 \pi m\sqrt{1 - 1/\sigma}}\,.
$$
Thus we obtain the following

\begin{Lemma}
\label{Lemma:hatvarphiexp(N/2)}
Let $N \in 2 \mathbb N$ and $\sigma \ge \frac{5}{4}$ be given, where $N_1 = \sigma N \in 2 \mathbb N$. Further let $m \in \mathbb N$ with $2m \ll N_1$, $\beta = b m$, and $b = 2\pi\big(1 - \frac{1}{2 \sigma}\big)$.

Then we have
\begin{eqnarray*}
{\hat \varphi}_{\mathrm{cexp}}\big(\frac{N}{2}\big)&\ge& \frac{2 m}{({\mathrm e}^{\beta} - 1)\,N_1}\Big[\frac{b}{5\, \sqrt{2 \pi m}}\,\big(1 - \frac{1}{\sigma}\big)^{-3/4} \,{\mathrm e}^{2 \pi m\sqrt{1 - 1/\sigma}} -\, 1 - \gamma(m,\sigma)\Big]\,.
\end{eqnarray*}
\end{Lemma}

By \eqref{eq:esigmaNvarphi}, the constant $e_{\sigma,N}({\hat \varphi}_{\mathrm{cexp}})$ can be estimated as follows
$$
e_{\sigma,N}({\hat \varphi}_{\mathrm{cexp}}) \le \frac{1}{{\hat \varphi}_{\mathrm{cexp}}(N/2)}\, \max_{n\in I_N}\big\| \sum_{r \in \mathbb Z\setminus \{0\}} {\hat \varphi}_{\mathrm{cexp}}(n + r N_1)\, {\mathrm e}^{2 \pi {\mathrm i}\,r N_1 \, \cdot}\big\|_{C(\mathbb T)}\,,
$$
where it holds by \eqref{eq:splitvarphicexp},
$$
\big\| \sum_{r \in \mathbb Z\setminus \{0\}} {\hat \varphi}_{\mathrm{cexp}}(n + r N_1)\, {\mathrm e}^{2 \pi {\mathrm i}\,r N_1 \, \cdot}\big\|_{C(\mathbb T)}
\le\, \sum_{r \in \mathbb Z\setminus \{0\}}|{\hat \psi}(n + r N_1)| +  \sum_{r \in \mathbb Z\setminus \{0\}} |{\hat \rho}(n + r N_1)|\,.
$$

\begin{Lemma}
\label{Lemma:sumvarphiexp}
Let $N \in 2 \mathbb N$ and $\sigma \ge \frac{5}{4}$ be given, where $N_1 = \sigma N \in 2 \mathbb N$. Further let $m \in \mathbb N$ with $2m \ll N_1$, $\beta = b m$, and $b = 2 \pi \big(1 - \frac{1}{2 \sigma}\big)$.

Then it holds for all $n\in I_N$,
$$
 \sum_{r\in \mathbb Z \setminus \{0\}} |{\hat \psi}(n + r N_1)|\le \, \frac{\beta}{({\mathrm e}^{\beta} - 1)\,N_1}\Big[2 \pi m + \frac{10}{\sqrt{2\pi m}}\big(1 - \frac{1}{2\sigma}\big)^{-1/2} \Big]\,.
$$
\end{Lemma}

\emph{Proof}. For all $n \in I_N$ and $r \in \mathbb Z \setminus \{0\}$ with $(n,r) \not= \big(-\frac{N}{2},1\big)$, it follows that
$$
| n + r N_1| > N_1 - \frac{N}{2} = N_1\,\big(1 - \frac{1}{2\sigma}\big)
$$
and hence by $\sigma \ge \frac{5}{4}$,
$$
\big| \frac{n}{N_1} + r \big| > 1 - \frac{1}{2 \sigma} \ge \frac{3}{5}\,.
$$
Thus by \eqref{eq:FTpsi}, we receive
$$
{\hat \psi}(n + r N_1) = \frac{\beta}{(e^{\beta} - 1)\,N_1}
\Big(\big(\frac{n}{N_1} + r\big)^2 - \big(1 - \frac{1}{2\sigma}\big)^2\Big)^{-1/2}\, J_1\big(2 \pi m \sqrt{\big(\frac{n}{N_1} + r\big)^2 - \big(1 - \frac{1}{2\sigma}\big)^2}\big)\,.
$$
In the case $(n,r) = \big(-\frac{N}{2},1\big)$, from \eqref{eq:FTpsi} it follows that
$$
{\hat \psi}\big(- \frac{N}{2} +  N_1\big) = {\hat \psi}\big(N_1 \big(1 - \frac{1}{2\sigma}\big) \big) = \frac{\pi m \beta}{2 ({\mathrm e}^{\beta} - 1)\,N_1}\,.
$$
Now we estimate the sum
\begin{eqnarray*}
\sum_{r \in \mathbb Z \setminus \{0\}} |{\hat \psi}(n + r N_1)| &=& |{\hat \psi}(n +  N_1)| + |{\hat \psi}(n - N_1)|\\
& & + \,\sum_{r \in \mathbb Z \setminus \{0,\pm 1\}} |{\hat \psi}(n + r N_1)|\,.
\end{eqnarray*}
For $n\in I_N$ and $r \in \mathbb Z \setminus \{0,\pm 1\}$ we have by  $\sigma \ge \frac{5}{4}$,
$$
\big| \frac{n}{N_1} + r \big| \ge 2 - \frac{1}{2 \sigma} \ge \frac{7}{5}
$$
and therefore by $m \ge 2$,
\begin{eqnarray*}
2 \pi m \sqrt{\big(\frac{n}{N_1} + r\big)^2 - \big(1 - \frac{1}{2\sigma}\big)^2} &\ge& 4 \pi \sqrt{\big(2 -  \frac{1}{2 \sigma}\big)^2 - \big(1 - \frac{1}{2 \sigma}\big)^2}\\
&=& 4 \pi \sqrt {3 - \frac{1}{\sigma}} \ge 4 \pi \sqrt{\frac{11}{5}} > 6\,.
\end{eqnarray*}
By an inequality for the Bessel function $J_1(x)$ (see \cite{Kra} or \cite[Lemma 3.2]{PoTa20}), it holds for all $x \ge 6$,
$$
|J_1(x)| \le \frac{1}{\sqrt x}
$$
such that for $x = 2 \pi m \sqrt{(r + n/N_1)^2 - (1 - 1/(2\sigma))^2}$,
\begin{eqnarray*}
& &\big|J_1\big(2 \pi m \sqrt{\big(\frac{n}{N_1} + r\big)^2 - \big(1 - \frac{1}{2\sigma}\big)^2}\big)\big|\\
& & \le\, \frac{1}{\sqrt{2 \pi m}}\,\Big(\big(\frac{n}{N_1} + r\big)^2 - \big(1 - \frac{1}{2 \sigma}\big)^2\Big)^{-1/4}\,.
\end{eqnarray*}
The decreasing function $h:\,\big[\frac{5}{4},\,\infty\big) \to \mathbb R$ defined by
$$
h(x) := x^{3/2}\,\big(x^2 - \big(1 - \frac{1}{2 \sigma}\big)^2\big)^{-3/4} = \Big(1 - \big(1 - \frac{1}{2 \sigma}\big)^2\,x^{-2}\Big)^{-3/4}\,,
$$
is bounded by
$$
h\big(\frac{5}{4}\big) = \Big(1 - \frac{16}{25}\,\big(1 - \frac{1}{2\sigma}\big)^2\Big)^{-3/4} \le \big(1 - \frac{16}{25}\big)^{-3/4} = \big(\frac{5}{3}\big)^{3/2} < \frac{5}{2}\,.
$$
Thus we receive for all $n \in I_N$ and $r \in \mathbb Z \setminus \{0, \pm 1\}$,
\begin{eqnarray*}
& &\Big(\big(\frac{n}{N_1} + r\big)^2 - \big(1 - \frac{1}{2\sigma}\big)^2 \Big)^{-1/2}\, \big|J_1\big(2 \pi m \sqrt{\big(\frac{n}{N_1} + r\big)^2 - \big(1 - \frac{1}{2\sigma}\big)^2}\big)\big|\\
& &\le \, \frac{1}{\sqrt{2 \pi m}}\Big(\big(\frac{n}{N_1} + r\big)^2 - \big(1 - \frac{1}{2\sigma}\big)^2\Big)^{-3/4}
\le \frac{5}{2\,\sqrt{2 \pi m}}\, \big|\frac{n}{N_1} + r \big|^{-3/2}\,.
\end{eqnarray*}
Hence for all $n\in I_N$, we obtain by \eqref{eq:sumxmu} that
\begin{eqnarray*}
\sum_{r \in \mathbb Z \setminus \{0,\pm 1\}} |{\hat \psi}(n + r N_1)| &\le& \frac{5 \beta}{2\,\sqrt{2 \pi m}\,({\mathrm e}^{\beta} - 1)\,N_1}\sum_{r \in \mathbb Z \setminus \{0,\pm 1\}} \big| \frac{n}{N_1} + r\big|^{-3/2}\\
&\le& \frac{10 \beta}{\sqrt{2 \pi m}\,({\mathrm e}^{\beta} - 1)\,N_1}\, \big(1 - \frac{1}{2\sigma}\big)^{-1/2}\,.
\end{eqnarray*}
Now we estimate ${\hat \psi}(v)$ for $v= n \pm N_1$, $n \in I_N$. For $v = - \frac{N}{2} + N_1 = N_1\big(1 - \frac{1}{2\sigma}\big)$, it holds by \eqref{eq:FTpsi},
\begin{equation}
\label{eq:FTpsispecial}
{\hat \psi}\big( - \frac{N}{2} + N_1\big) = {\hat \psi}\big( N_1\,\big(1 - \frac{1}{2 \sigma}\big)\big) = \frac{\pi m \beta}{2\,({\mathrm e}^{\beta} - 1)\,N_1}\,.
\end{equation}
For all the other $v= n \pm N_1 \not= - \frac{N}{2} + N_1$, $n \in I_N$, we have
$$
|n \pm N_1| > \big(1 - \frac{1}{2\sigma}\big)\,N_1
$$
such that by \eqref{eq:FTpsi}, ${\hat \psi}(n \pm N_1)$ reads as follows
$$
\frac{\beta}{({\mathrm e}^{\beta}-1)\,N_1} \Big(\big(1 \pm \frac{n}{N_1}\big)^2 - \big(1 - \frac{1}{2\sigma}\big)^2\Big)^{-1/2}\,J_1\big(2\pi m \sqrt{\big(1 \pm \frac{n}{N_1}\big)^2 - \big(1 - \frac{1}{2\sigma}\big)^2}\big)\,.
$$
Since
$$
2\pi m \sqrt{\big(1 \pm \frac{n}{N_1}\big)^2 - \big(1 - \frac{1}{2\sigma}\big)^2}  > 0
$$
can be small for $n \in I_N$, we estimate the Bessel function $J_1(x)$, $x\ge 0$, by Poisson's integral (see \cite[p.~47]{Watson})
\begin{equation*}
|J_1(x)| = \frac{x}{\pi}\,\big|\int_0^{\pi} \cos(x\, \cos t)\,(\sin t)^2\,{\mathrm d}t\big| \le  \frac{x}{\pi}\,\int_0^{\pi} (\sin t)^2\,{\mathrm d}t =\frac{x}{2}\,,
\end{equation*}
such that
$$
|{\hat \psi}(n \pm N_1)| \le \frac{\pi m \beta}{({\mathrm e}^{\beta} - 1)\,N_1}\,.
$$
By \eqref{eq:FTpsispecial}, this estimate of $|{\hat \psi}(n \pm N_1)|$ is valid for all $n \in I_N$.
This completes the proof. \qedsymbol
\medskip

Now for arbitrary $n \in I_N$, we have to estimate the series
$$
\sum_{r \in \mathbb Z\setminus \{0\}} |{\hat \rho}(n + r N_1)|\,.
$$
By \eqref{eq:hatrhoinw} and Lemma \ref{Lemma:est|I(w)|}, we obtain for any $v \in \mathbb R \setminus \{0\}$,
\begin{eqnarray*}
|{\hat \rho}(v)| &\le& \frac{2 m}{({\mathrm e}^{\beta}-1)\,N_1}\Big[  \big(1 + \frac{1}{\mathrm e}\big)\,\frac{\beta}{2\pi m}\,\sqrt[4]{\frac{5}{2\pi m}}\,\big(\frac{|v|}{N_1}\big)^{-5/4}\\
 & & \,+ \,\frac{N_1}{\pi m\,|v|}\,{\mathrm e}^{-\sqrt{2 \pi m\,|v|/N_1}} + \big({\mathrm e}^{-\sqrt 2 \beta} + \frac{1}{2}\big)\,{\mathrm e}^{-2\pi m\,|v|/N_1}\Big]\,.
\end{eqnarray*}
Thus we obtain that
\begin{eqnarray*}
\sum_{r \in \mathbb Z\setminus \{0\}} |{\hat \rho}(n + r N_1)|&\le& \frac{2 m}{({\mathrm e}^{\beta}-1)\,N_1}\Big[  \big(1 + \frac{1}{\mathrm e}\big)\,\frac{\beta}{2\pi m}\,\sqrt[4]{\frac{5}{2\pi m}}\,S_1(n)\\
& & +\, 2\,S_2(n) + \big({\mathrm e}^{-\sqrt 2 \beta} + \frac{1}{2}\big)\,S_3(n)\Big]
\end{eqnarray*}
with
\begin{eqnarray}
S_1(n) &:=& \sum_{r \in \mathbb Z\setminus \{0\}} \big|\frac{n}{N_1}+r\big|^{-5/4}\,,\label{eq:S1(n)}\\
S_2(n) &:=& \sum_{r \in \mathbb Z\setminus \{0\}} \frac{1}{2 \pi m\,|r + n/N_1|}\,{\mathrm e}^{-\sqrt{2 \pi m\,|r + n/N_1|}}\,,\label{eq:S2(n)}\\
S_3(n) &:=& \sum_{r \in \mathbb Z\setminus \{0\}} {\mathrm e}^{-2 \pi m\,|r + n/N_1|}\,. \label{eq:S3(n)}
\end{eqnarray}
The inequalities \eqref{eq:sumxmu}, \eqref{eq:sumexp(-ax)}, and \eqref{eq:sumexp(-sqrtx)/x}      imply that
\begin{eqnarray}
S_1(n) &\le& \big(\frac{4\sigma}{2\sigma -1} + 8\big)\, \big(1 - \frac{1}{2 \sigma}\big)^{-1/4}\,, \label{eq:estS1(n)}\\
S_2(n) &\le& \frac{2 \sigma}{(2\sigma -1)\,\pi m}\,{\mathrm e}^{-\sqrt{2 \pi m - \pi m/\sigma}} + \frac{2}{\pi m}\,E_1\big(\sqrt{2 \pi m - \pi m/\sigma}\big)\,, \label{eq:estS2(n)}\\
S_3(n) &\le& \big(2 + \frac{1}{\pi m}\big)\,{\mathrm e}^{-2 \pi m + \pi m/\sigma}\,. \label{eq:estS3(n)}
\end{eqnarray}
Thus we obtain the following

\begin{Lemma}
\label{Lemma:estsumhatrho}
Let $N\in 2 \mathbb N$ and $\sigma \ge \frac{5}{4}$ be given, where $N_1 = \sigma\,N \in 2 \mathbb N$. Further let $m\in \mathbb N\setminus \{1\}$ with $2 m \ll N_1$,  $\beta = b m$, and $b = 2 \pi \big(1 - \frac{1}{2\sigma}\big)$.

Then for each $n \in I_N$, it holds the estimate
\begin{eqnarray*}
\sum_{r\in \mathbb Z \setminus \{0\}} |{\hat \rho}(n + rN_1)| &\le& \frac{1}{({\mathrm e}^{\beta} -1)\,N_1}\Big[\big(1 + \frac{1}{\mathrm e}\big)\,\frac{\beta}{\pi }\,\sqrt[4]{\frac{5}{2\pi m}}\,\big(\frac{4\sigma}{2\sigma -1} + 8\big)\, \big(1 - \frac{1}{2 \sigma}\big)^{-1/4}\\
& & + \,\frac{8 \sigma}{(2\sigma -1)\,\pi}\,{\mathrm e}^{-\sqrt{2 \pi m - \pi m/\sigma}} + \frac{8}{\pi}\,E_1\big(\sqrt{2 \pi m - \pi m/\sigma}\big)\\
& &+ \,\big({\mathrm e}^{-\sqrt 2 \beta} + \frac{1}{2}\big)\, \big(4m + \frac{2}{\pi}\big)\,{\mathrm e}^{-2 \pi m + \pi m/\sigma}\Big]\,.
\end{eqnarray*}
\end{Lemma}

Hence from Lemmas \ref{Lemma:sumvarphiexp} and \ref{Lemma:estsumhatrho} it follows that
\begin{equation}
\label{eq:normsumhatvarphicexp}
\max_{n \in I_N} \,\big\| \sum_{r \in {\mathbb Z}\setminus \{0\}} {\hat \varphi}_{\mathrm{cexp}}(n +rN_1)\,{\mathrm e}^{2 \pi {\mathrm i} r N_1\, \cdot}\big\|_{C(\mathbb T)} \le \frac{\beta}{({\mathrm e}^{\beta} - 1)\,N_1}\, b(m,\sigma)
\end{equation}
with the constant
\begin{eqnarray}
\label{eq:b(msigma)}
b(m,\sigma) &:=& 2 \pi m + \frac{10}{\sqrt{2\pi m}}\,\big(1 - \frac{1}{2\sigma}\big)^{-1/2} \nonumber \\
& & + \,\big(1 + \frac{1}{\mathrm e}\big)\,\frac{1}{\pi}\,\sqrt[4]{\frac{5}{2\pi m}}\,\big(\frac{4\sigma}{2\sigma -1} + 8\big)\, \big(1 - \frac{1}{2 \sigma}\big)^{-1/4} \nonumber\\
& & + \,\frac{8 \sigma}{(2\sigma -1)\,\pi \beta}\,{\mathrm e}^{-\sqrt{2 \pi m - \pi m/\sigma}} + \frac{8}{\pi \beta}\,E_1\big(\sqrt{2 \pi m - \pi m/\sigma}\big) \nonumber\\
& & + \,\frac{4 \pi m + 2}{\beta \pi}\,\big({\mathrm e}^{-\sqrt 2 \beta} + \frac{1}{2}\big)\,{\mathrm e}^{-2 \pi m + \pi m/\sigma}\,.
\end{eqnarray}
Using Lemma \ref{Lemma:hatvarphiexp(N/2)}, we obtain by
$$
e_{\sigma,N}(\varphi_{\mathrm{cexp}}) \le \frac{1}{{\hat \varphi}_{\mathrm{cexp}}(N/2)}\,\max_{n \in I_N} \big\| \sum_{r \in {\mathbb Z}\setminus \{0\}} {\hat \varphi}_{\mathrm{cexp}}(n +rN_1)\,{\mathrm e}^{2 \pi {\mathrm i} r N_1\, \cdot}\big\|_{C(\mathbb T)}
$$
the following

\begin{Theorem}
\label{Theorem:expwindow}
Let $N \in 2 \mathbb N$ and $\sigma \ge \frac{5}{4}$ be given, where $N_1 = \sigma N \in 2 \mathbb N$. Further let $m \in \mathbb N$ with $2m \ll N_1$, $\beta = b m$, and $b = 2\pi \big(1 - \frac{1}{2\sigma}\big)$.

Then the $C(\mathbb T)$-error constant of the continuous $\exp$-type window function \eqref{eq:expwindow} can be estimated by
$$
e_{\sigma}(\varphi_{\mathrm{cexp}}) \le \frac{\beta\,b(m,\sigma)}{2m}\,\Big[\frac{b}{5\, \sqrt{2 \pi m}} \,\big(1 - \frac{1}{\sigma}\big)^{-3/4}\,{\mathrm e}^{2 \pi m\sqrt{1 - 1/\sigma}} - 1 - \gamma(m,\sigma)\Big]^{-1}
$$
In other words, the continuous $\exp$-type window function \eqref{eq:expwindow} is convenient for $\mathrm{NFFT}$.
\end{Theorem}

Note that for $\sigma \in \big[\frac{5}{4},\,2\big]$ and $m \ge 2$, it holds by \eqref{eq:b(msigma)},
$$
\frac{\beta\,b(m,\sigma)}{2m} = \pi b(m,\sigma)\,\big(1 - \frac{1}{2\sigma}\big) \le \frac{3 \pi^2}{2} \,m + b_0 < 15\,m + b_0
$$
with $b_0<17$.

In order to compute the Fourier transform $\hat \varphi$ of  window function $\varphi \in \Phi_{m,N_1}$, we  approximate this window function by numerical integration. In our next numerical examples we apply the following method.
Since the window function $\varphi\in \Phi_{m,N_1}$ is even and supported in
$[-\frac{m}{N_1}, \frac{m}{N_1}]$, we have
\begin{eqnarray*}
	{\hat \varphi}(v) &=& \int_{\mathbb R} \varphi(x)\,{\mathrm e}^{-2\pi {\mathrm i}\,v x}\, {\mathrm d}x
	= \frac{m}{N_1}\,\int_{-1}^1
	\varphi\big(\frac{m}{N_1}\, t\big)\,
	{\mathrm e}^{-2\pi {\mathrm i}\,m v t/N_1}\, {\mathrm d}t\\ [1ex]
	&=& \frac{2m}{N_1}\,\int_{0}^1
	\varphi\big(\frac{m}{N_1}\, t\big)\,
	\cos\frac{2\pi \,m v t}{N_1}\, {\mathrm d}t\,.
\end{eqnarray*}

We evaluate the last integral using a global adaptive quadrature \cite{Sham08} for $\hat \varphi(k)$, $k=0,\ldots, N$. In general, this values can be precomputed, see \cite{nfft3, FINUFFT}.

\begin{Remark}
\label{Remark:origexpwindow}
As in \cite{BaMaKl18, Ba20}, we can consider also the original $\exp$\emph{-type window function}
\begin{equation}
\label{eq:origexpwindow}
\varphi_{\exp}(x) := \left\{ \begin{array}{ll} \exp\big(\beta \sqrt{1 - (N_1x/m)^2} - \beta\big) & \quad x\in I\,,\\ [1ex]
\frac{1}{2}\,{\mathrm e}^{-\beta} & \quad x = \pm \frac{m}{N_1}\,,\\ [1ex]
0 & \quad x \in \mathbb R \setminus \bar I
\end{array} \right.
\end{equation}
with the shape parameter $\beta = m b = 2\pi m \big(1 - \frac{1}{2\sigma}\big)$, $\sigma \ge \frac{5}{4}$, $N \in 2 \mathbb N$, and $N_1 = \sigma N \in 2 \mathbb N$. Further we assume that $m \in \mathbb N \setminus \{1\}$ fulfills $2 m \ll N_1$. This window function possesses jump discontinuities at $x = \pm \frac{m}{N_1}$ with very small jump height ${\mathrm e}^{-\beta}$, such that \eqref{eq:origexpwindow} is ``almost continuous''.

We split \eqref{eq:origexpwindow} in the form
$$
\varphi_{\exp}(x) = \big(1 - {\mathrm e}^{-\beta}\big)\,\varphi_{\mathrm{cexp}}(x) + {\mathrm e}^{-\beta}\,\varphi_{\mathrm{rect}}(x)\,, \quad x \in \mathbb R\,,
$$
with the window functions \eqref{eq:expwindow} and \eqref{eq:rectwindow}. Then the Fourier transform of \eqref{eq:origexpwindow} reads as follows
$$
{\hat \varphi}_{\exp}(v) = \big(1 - {\mathrm e}^{-\beta}\big)\,{\hat \varphi}_{\mathrm{cexp}}(v) + {\mathrm e}^{-\beta}\,{\hat \varphi}_{\mathrm{rect}}(v)\,, \quad v \in \mathbb R\,,
$$
By Lemma \ref{Lemma:hatvarphiexp(N/2)} it follows that
\begin{eqnarray*}
{\hat \varphi}_{\mathrm{exp}}\big(\frac{N}{2}\big) &=& \min_{n\in I_N} {\hat \varphi}_{\mathrm{exp}}(n)\\
&=&  \min_{n\in I_N} \Big[\big(1 - {\mathrm e}^{-\beta}\big)\,{\hat \varphi}_{\mathrm{cexp}}(n) + {\mathrm e}^{-\beta}\,{\hat \varphi}_{\mathrm{rect}}(n) \Big]\\
&\ge& \frac{2 m}{N_1 \,{\mathrm e}^{\beta}}\Big[\frac{b}{5\, \sqrt{2 \pi m}}\, \big(1 -\frac{1}{\sigma}\big)^{-3/4}\,{\mathrm e}^{2 \pi m\sqrt{1 - 1/\sigma}} - \gamma(m,\sigma)\Big]\,.
\end{eqnarray*}
Using \eqref{eq:normsumhatvarphicexp} and \eqref{eq:sumhatphirect1}, we estimate for all $n \in I_N$,
\begin{eqnarray*}
& &\big\|\sum_{r \in \mathbb Z \setminus \{0\}} {\hat \varphi}_{\exp}(n + r N_1)\,{\mathrm e}^{2\pi {\mathrm i}r N_1 \, \cdot}\big\|_{C(\mathbb T)} \le \big(1 - {\mathrm e}^{-\beta}\big)\sum_{r \in \mathbb Z \setminus \{0\}} |{\hat \varphi}_{\mathrm{cexp}}(n + r N_1)|\\
& & + \, {\mathrm e}^{-\beta}\,\big\|\sum_{r \in \mathbb Z \setminus \{0\}} {\hat \varphi}_{\mathrm{rect}}(n + r N_1)\,{\mathrm e}^{2\pi {\mathrm i}r N_1 \, \cdot}\big\|_{C(\mathbb T)}\\
&  &\le\, \frac{1}{N_1\,{\mathrm e}^{\beta}} \,(\beta\,b(m,\sigma) + 3m)\,.
\end{eqnarray*}
Thus we obtain
$$
e_{\sigma}(\varphi_{\exp}) \le \big(\frac{\beta\, b(m,\sigma)}{2 m} + \frac{3}{2}\big)\,\Big[\frac{b}{5 \sqrt{2 \pi m} } \,\big(1 - \frac{1}{\sigma}\big)^{-3/4}\,{\mathrm e}^{2 \pi m\sqrt{1 - 1/\sigma}} -\, \gamma(m,\sigma)\Big]^{-1}\,.
$$
Thus the discontinuous window function \eqref{eq:origexpwindow} possesses a similar $C(\mathbb T)$-error constant as the continuous $\exp$-type window function \eqref{eq:expwindow}.
\end{Remark}

\section{Continuous cosh-type window function}\label{Sec:coshWindow}

For fixed shape parameter $\beta = bm = 2\pi m\big(1- \frac{1}{2\sigma}\big)$ and oversampling factor $\sigma \ge \frac{5}{4}$, we consider the \emph{continuous} $\cosh$-\emph{type window function}
\begin{equation}
\label{eq:coshwindow}
\varphi_{\cosh}(x) := \left\{ \begin{array}{ll}
\frac{1}{\cosh \beta - 1}\,\Big(\cosh\left(\beta \sqrt{1- \left(\frac{N_1 x}{m}\right)^2 }\right) - 1\Big) & \quad x \in I\,, \\ [1ex]
0 & \quad x \in \mathbb R \setminus  I\,.
\end{array}\right.
\end{equation}
Obviously, $\varphi_{\cosh} \in \Phi_{m,N_1}$ is a continuous window function. Note that recently a discontinuous version of this window function was suggested in \cite[Remark 13]{BaMaKl18}. But up to now, a corresponding error estimate for the related NFFT was unknown. Now we show that the $C(\mathbb T)$-error constant $e_{\sigma}(\varphi_{\cosh})$ can be estimated by a similar upper bound as $e_{\sigma}(\varphi_{\mathrm{cexp}})$ in Theorem \ref{Theorem:expwindow}. Thus
the window functions \eqref{eq:expwindow} and \eqref{eq:coshwindow} possess the same error behavior with respect to the NFFT.
\medskip

In the following, we use the same technique as in Section \ref{Sec:expWindow}.
Since the Fourier transform of \eqref{eq:coshwindow} is unknown analytically, we represent \eqref{eq:coshwindow} as the sum
$$
\varphi_{\cosh}(x) = \psi_1(x) + \rho_1(x)\,,
$$
where the Fourier transform of $\psi_1$ is known and where the correction term $\rho_1$ has small magnitude $|\rho_1|$. We choose
$$
\psi_1(x) := \left\{ \begin{array}{ll}
\frac{1}{\cosh \beta - 1} \, \sinh \big(\beta \,\sqrt{1 - (N_1 x)^2/m^2}\big) & \quad x \in I\,, \\
0 & \quad x \in \mathbb R \setminus  I
\end{array}\right.
$$
and
$$
\rho_1(x) := \left\{ \begin{array}{ll}
\frac{1}{\cosh \beta - 1} \, \Big({\mathrm e}^{-\beta \,\sqrt{1 - (N_1 x)^2/m^2}}-1\Big) & \quad x \in I\,, \\ [1ex]
0 & \quad x \in \mathbb R \setminus  I\,.
\end{array}\right.
$$
Since $\rho_1$ is even and $\rho_1 |_{ [0,\, \frac{m}{N_1}]}$ is increasing, we have
$$
0 \ge \rho_1(x) \ge \rho_1(0) = \frac{{\mathrm e}^{-\beta}-1}{\cosh \beta - 1} = - \frac{2}{{\mathrm e}^{\beta}-1}\,.
$$
Since $\rho_1$ has a small values in the compact support $\bar I$, the Fourier transform $\hat \rho_1$ is small too and it holds
$$
|{\hat \rho_1}(v)| = \big|\int_I \rho(x)\,{\mathrm e}^{2\pi {\mathrm i}\,v x}\,{\mathrm d}x \big| \le \frac{4 m}{({\mathrm e}^{\beta} - 1)\,N_1}\,.
$$
Substituting $t= N_1 x/m$, we determine the Fourier transform
\begin{eqnarray*}
{\hat \varphi}_{\cosh}(v) &=& \int_I \varphi_{\cosh}(x)\,{\mathrm e}^{-2\pi {\mathrm i}\,v x}\, {\mathrm d}x
= \int_I \psi_1(x)\,{\mathrm e}^{-2\pi {\mathrm i}\,v x}\, {\mathrm d}x + \int_I \rho_1(x)\,{\mathrm e}^{-2\pi {\mathrm i}\,v x}\, {\mathrm d}x \\
&=& \frac{2 m}{(\cosh \beta - 1)\,N_1}\Big[\int_0^1 \sinh\big(\beta \sqrt{1-t^2}\big)\,\cos\frac{2\pi m v t}{N_1}\, {\mathrm d}t\\
& &+\, \int_0^1 \big({\mathrm e}^{-\beta \sqrt{1 - t^2}}-1\big) \, \cos\frac{2\pi m v t}{N_1}\,{\mathrm d}t\Big]\,.
\end{eqnarray*}
For simplicity, we introduce the scaled frequency  $w:= 2 \pi m v/N_1$ such that
\begin{eqnarray*}
{\hat \varphi}_{\exp}\big(\frac{N_1 w}{2 \pi m}\big) &=& \frac{2 m}{(\cosh \beta - 1)\,N_1} \Big[\int_0^1 \sinh\big(\beta \sqrt{1-t^2}\big)\,\cos(w t)\, {\mathrm d}t\nonumber \\
& &+\, \int_0^1 \big({\mathrm e}^{-\beta \sqrt{1 - t^2}}- 1\big) \, \cos(w t)\,{\mathrm d}t\Big]\,.
\end{eqnarray*}
From \cite[p. 38]{Oberh90} it follows that
$$
\int_0^1 \sinh\big(\beta \sqrt{1-t^2}\big)\,\cos(w t)\, {\mathrm d}t = \frac{\pi \beta}{2} \left\{\begin{array}{ll}
(\beta^2 - w^2)^{-1/2}\,I_1\big(\sqrt{\beta^2 - w^2}\big) & \quad w \in (-\beta,\, \beta)\,,\\
1/4 & \quad w = \pm \beta\,, \\
(w^2 - \beta^2)^{-1/2}\,J_1\big(\sqrt{w^2 - \beta^2}\big) & \quad w \in \mathbb R \setminus [-\beta,\, \beta]\,,
\end{array} \right.
$$
where $I_1$ denotes the modified Bessel function and $J_1$ the Bessel function of first order.

Using Lemma \ref{Lemma:estint01bysinc}, we receive for $w \in (-\beta,\,\beta)$,
$$
{\hat \varphi}_{\cosh}\big(\frac{N_1 w}{2 \pi m}\big)\ge \frac{2 m}{(\cosh \beta - 1)\,N_1}\Big[\frac{\pi \beta}{2}\,(\beta^2 - w^2)^{-1/2}\,I_1\big(\sqrt{\beta^2 - w^2}\big) - \mathrm{sinc}\,w - \gamma(m,\sigma)\Big]\,.
$$
This means for $\beta = b m$ and $v = \frac{N_1 w}{2 \pi m}\in \big(-N_1\, \big(1 - \frac{1}{2\sigma}\big),\,N_1\, \big(1 - \frac{1}{2\sigma}\big)\big)$,
\begin{eqnarray*}
{\hat \varphi}_{\cosh}(v) &\ge& \frac{2 m}{(\cosh \beta - 1)\,N_1}\Big[ \frac{b}{4}\,\big(\big(1- \frac{1}{2\sigma}\big)^2 - \frac{v^2}{N_1^2}\big)^{-1/2}\,I_1\big(2\pi m \sqrt{\big(1-\frac{1}{2\sigma}\big)^2 - \frac{v^2}{N_1^2}}\big)\\
& & - \, \mathrm{sinc}\frac{2 \pi m v}{N_1} - \gamma(m,\sigma)\Big]\,.
\end{eqnarray*}
Since the function $h:\,[0,\beta) \to \mathbb R$ defined by
$$
h(w) := \frac{\pi \beta}{2}\,\big(\beta^2 - w^2\big)^{-1/2}\,I_1\big(\sqrt{\beta^2 - w^2}\big) - \mathrm{sinc}\,w > 0
$$
is decreasing and
$$
\lim_{w \to \beta-0} h(w) = \frac{\pi \beta}{4} - \mathrm{sinc}\, \beta \ge \pi^2\,\big(1 - \frac{1}{2 \sigma}\big) - 1 \ge \frac{3\pi^2}{5}- 1\,,
$$
the Fourier transform ${\hat \varphi}_{\cosh}(v)$ is positive and decreasing for $v \in \big[0,\, N_1\,\big(1 - \frac{1}{2 \sigma}\big)\big)$ too. Hence we obtain
\begin{eqnarray*}
\min_{n\in I_N} {\hat \varphi}_{\cosh}(n) &=& {\hat \varphi}_{\cosh}\big(\frac{N}{2}\big)\\
&\ge& \frac{2 m}{(\cosh \beta - 1)\,N_1}\Big[\frac{b}{4}\,\big(1 - \frac{1}{\sigma}\big)^{-1/2}\,I_1\big(2 \pi m \sqrt{1 - \frac{1}{\sigma}}\big) - \mathrm{sinc}\frac{\pi m}{\sigma} - \gamma(m,\sigma)\Big]\\
&\ge& \frac{2 m}{(\cosh \beta - 1)\,N_1}\Big[\frac{b}{4}\,\big(1 - \frac{1}{\sigma}\big)^{-1/2}\,I_1\big(2 \pi m \sqrt{1 - \frac{1}{\sigma}}\big) - 1 - \gamma(m,\sigma)\Big]\,.
\end{eqnarray*}
From $m \ge 2$ and $\sigma \ge \frac{5}{4}$ it follows that
$$
2 \pi m \sqrt{1 - \frac{1}{\sigma}} \ge 4 \pi\,\sqrt{1 - \frac{1}{\sigma}} \ge x_0 := \frac{4 \pi}{\sqrt 5}\,.
$$
By the inequality (see \cite{Bar} or \cite[Lemma 3.3]{PoTa20})
$$
\sqrt{2 \pi x_0}\,{\mathrm e}^{-x_0}\, I_1(x_0) \le \sqrt{2 \pi x}\,{\mathrm e}^{-x}\, I_1(x)\,, \quad x \ge x_0\,,
$$
we sustain that
$$
I_1(x) \ge \sqrt x_0 \,{\mathrm e}^{-x_0}\,I_1(x_0)\,x^{-1/2}\,{\mathrm e}^{x} > \frac{2}{5}\,x^{-1/2}\,{\mathrm e}^{x}\,, \quad x \ge x_0\,.
$$
Hence for $x = 2 \pi m\,\sqrt{1 - 1/\sigma}$ we get the estimate
$$
I_1\big(2 \pi m \sqrt{1 - \frac{1}{\sigma}}\big) > \frac{\sqrt 2}{5\, \sqrt{\pi m}} \,\big(1 -\frac{1}{\sigma}\big)^{-1/4}\,{\mathrm e}^{2 \pi m\sqrt{1 - 1/\sigma}}\,.
$$
Thus we obtain the following

\begin{Lemma}
\label{Lemma:hatvarphicosh(N/2)}
Let $N \in 2 \mathbb N$ and $\sigma \ge \frac{5}{4}$ be given, where $N_1 = \sigma N \in 2 \mathbb N$. Further let $m \in \mathbb N$ with $2m \ll N_1$,  $\beta = b m$, and $b = 2 \pi \big(1 - \frac{1}{2 \sigma}\big)$.

Then we have
$$
{\hat \varphi}_{\cosh}\big(\frac{N}{2}\big)\le \frac{2m}{(\cosh \beta - 1)\,N_1}\Big[\frac{\sqrt \pi}{5\, \sqrt{2 m} }\,\big(1 - \frac{1}{2 \sigma}\big)\big(1 -\frac{1}{\sigma}\big)^{-3/4}\,{\mathrm e}^{2 \pi m\sqrt{1 - 1/\sigma}}- 1 - \gamma(m,\sigma)\Big]\,.
$$
\end{Lemma}

By \eqref{eq:esigmaNvarphi} the constant $e_{\sigma,N}({\hat \varphi}_{\cosh})$ can be estimated as follows
\begin{eqnarray*}
e_{\sigma,N}({\hat \varphi}_{\cosh}) &\le& \frac{1}{\min_{n\in I_N}{\hat \varphi}_{\cosh}(n)}\, \max_{n\in I_N}\big\| \sum_{r \in \mathbb Z\setminus \{0\}} {\hat \varphi}_{\cosh}(n + r N_1)\, {\mathrm e}^{2 \pi {\mathrm i}\,r N_1 \, \cdot}\big\|_{C(\mathbb T)}\\
& &=\, \frac{1}{{\hat \varphi}_{\cosh}(N/2)}\, \max_{n\in I_N}\big\| \sum_{r \in \mathbb Z\setminus \{0\}} {\hat \varphi}_{\cosh}(n + r N_1)\, {\mathrm e}^{2 \pi {\mathrm i}\,r N_1 \, \cdot}\big\|_{C(\mathbb T)}\,,
\end{eqnarray*}
where it holds
\begin{eqnarray}
& &\big\| \sum_{r \in \mathbb Z\setminus \{0\}} {\hat \varphi}_{\cosh}(n + r N_1)\, {\mathrm e}^{2 \pi {\mathrm i}\,r N_1 \, \cdot}\big\|_{C(\mathbb T)} \nonumber\\
& &\le\, \sum_{r \in \mathbb Z\setminus \{0\}}|{\hat \psi}_1(n + r N_1)| +  \sum_{r \in \mathbb Z\setminus \{0\}} |{\hat \rho}_1(n + r N_1)|\,. \label{eq:normsumhatphicosh}
\end{eqnarray}
Analogously to Lemma \ref{Lemma:sumvarphiexp}, we get

\begin{Lemma}
\label{Lemma:sumvarphicosh}
Let $N \in 2 \mathbb N$ and $\sigma \ge \frac{5}{4}$ be given, where $N_1 = \sigma N \in 2 \mathbb N$. Further let $m \in \mathbb N$ with $2m \ll N_1$, $\beta = b m$, and $b = 2 \pi \big(1 - \frac{1}{2 \pi}\big)$.

Then it holds for all $n \in I_N$,
$$
\sum_{r \in \mathbb Z\setminus \{0\}}|{\hat \psi}_1(n + r N_1)| \le \frac{\beta}{(\cosh \beta - 1)\,N_1}\,\Big[2 \pi m + \frac{10}{\sqrt{2 \pi m}}\,\big(1 - \frac{1}{2 \sigma}\big)^{-1/2}\Big]\,.
$$
\end{Lemma}

Finally for any $n \in I_N$, we estimate the sum
$$
\sum_{r \in \mathbb Z\setminus \{0\}}|{\hat \rho}_1(n + r N_1)|\,.
$$
Analogously to Lemma \ref{Lemma:estsumhatrho}, we obtain

\begin{Lemma}
\label{Lemma:sumhatrho1}
Let $N \in 2 \mathbb N$ and $\sigma \ge \frac{5}{4}$ be given, where $N_1 = \sigma N \in 2 \mathbb N$. Further let $m \in \mathbb N$ with $2m \ll N_1$, $\beta = b m$, and $b = 2 \pi \big(1 - \frac{1}{2 \sigma}\big)$.

Then for all $n \in I_N$,
\begin{eqnarray*}
\sum_{r \in \mathbb Z\setminus \{0\}}|{\hat \rho}_1(n + r N_1)| &\le& \frac{1}{\big(\cosh \beta - 1\big)\,N_1}\,\Big[\big(1 +\frac{1}{\mathrm e}\big) \frac{\beta}{\pi} \,\sqrt[4]{\frac{5}{2\pi m}}\,\big(\frac{4\sigma}{2\sigma -1} + 8\big)\, \big(1 - \frac{1}{2 \sigma}\big)^{-1/4}\\
& & + \,\frac{8 \sigma}{(2\sigma -1)\,\pi}\,{\mathrm e}^{-\sqrt{2 \pi m - \pi m/\sigma}} + \frac{8}{\pi }\,E_1\big(\sqrt{2 \pi m - \pi m/\sigma}\big)\\
& & + \,\big({\mathrm e}^{-\sqrt 2 \beta} + \frac{1}{2}\big)\, \big(4 m + \frac{2}{\pi}\big)\,{\mathrm e}^{-2 \pi m + \pi m/\sigma}\Big]\,.
\end{eqnarray*}
\end{Lemma}

From Lemmas  \ref{Lemma:sumvarphicosh} and \ref{Lemma:sumhatrho1} it follows that by \eqref{eq:normsumhatphicosh},
\begin{equation}
\label{eq:normsumhatphicosh1}
\big\| \sum_{r \in \mathbb Z\setminus \{0\}} {\hat \varphi}_{\cosh}(n + r N_1)\, {\mathrm e}^{2 \pi {\mathrm i}\,r N_1 \, \cdot}\big\|_{C(\mathbb T)} \le \frac{\beta\,b(m,\sigma)}{\big(\cosh \beta - 1\big)\,N_1}
\end{equation}
with the constant \eqref{eq:b(msigma)}.

Using
$$
e_{\sigma,N}(\varphi_{\cosh}) \le \frac{1}{{\hat \varphi}_{\cosh}(N/2)}\,\max_{n \in I_N}\big\| \sum_{r \in \mathbb Z\setminus \{0\}} {\hat \varphi}_{\cosh}(n + r N_1)\, {\mathrm e}^{2 \pi {\mathrm i}\,r N_1 \, \cdot}\big\|_{C(\mathbb T)}\,,
$$
it follows from \eqref{eq:normsumhatphicosh1} and  Lemma \ref{Lemma:hatvarphicosh(N/2)}:

\begin{Theorem}
\label{Theorem:coshwindow}
Let $N \in 2 \mathbb N$ and $\sigma \ge \frac{5}{4}$ be given, where $N_1 = \sigma N \in 2 \mathbb N$. Further let $m \in \mathbb N$ with $2m \ll N_1$, $\beta = b m$, and $b = 2 \pi \big(1 - \frac{1}{2 \sigma}\big)$.

Then the $C(\mathbb T)$-error constant of the continuous $\cosh$-type window function \eqref{eq:coshwindow} can be estimated by
$$
e_{\sigma}(\varphi_{\cosh}) \le \frac{\beta\,b(m,\sigma)}{2m}\,\Big[\frac{\sqrt \pi}{5\, \sqrt{2 m} }\, \big(1 - \frac{1}{2 \sigma}\big)\,\big(1 - \frac{1}{\sigma}\big)^{-3/4}\,{\mathrm e}^{2 \pi m\sqrt{1 - 1/\sigma}} - 1 - \gamma(m,\sigma)\Big]^{-1}\,,
$$
i.e., the continuous $\cosh$-type window function \eqref{eq:expwindow} is convenient for $\mathrm{NFFT}$.
\end{Theorem}

\section*{Conclusion}

In this paper, we prefer the use of continuous, compactly supported window functions for NFFT (with nonequispaced spatial data and equispaced frequencies). Such window functions simplify the algorithm for NFFT,
since the truncation error of NFFT vanishes. Further, such window functions can produce very small errors of NFFT. Examples of such window functions are the continuous Kaiser-Bessel window function \eqref{eq:Kaiser-Bessel},
continuous $\exp$-type window function \eqref{eq:expwindow}, $\sinh$-type window function \eqref{eq:sinhwindow}, and continuous $\cosh$-type window function \eqref{eq:coshwindow} which possess
 the same support and shape parameter. For these window functions, we present novel explicit error estimates for NFFT and
we derive rules for the convenient choice of the truncation parameter $m \ge 2$ and the oversampling parameter $\sigma \ge \frac{5}{4}$. The main tool of this approach is
the decay of the Fourier transform ${\hat \varphi}(v)$ of $\varphi \in \Phi_{m,N_1}$ for $|v| \to \infty$. A rapid decay of $\hat \varphi$ is essential for small error constants. Unfortunately,
the Fourier transform of certain window function  $\varphi$, such as \eqref{eq:expwindow} and \eqref{eq:coshwindow}, is unknown analytically. Therefore we propose a new technique and split $\varphi$ into a sum of two compactly supported functions $\psi$ and $\rho$, where the Fourier transform $\hat \psi$ is explicitly
known and where $|\rho|$ is sufficiently small. Further, it is shown that the standard Kaiser-Bessel window function and original $\exp$-type window function which have jump discontinuities with very small jump heights at the endpoints of their support, possess a similar error behavior as the corresponding  continuous window functions.

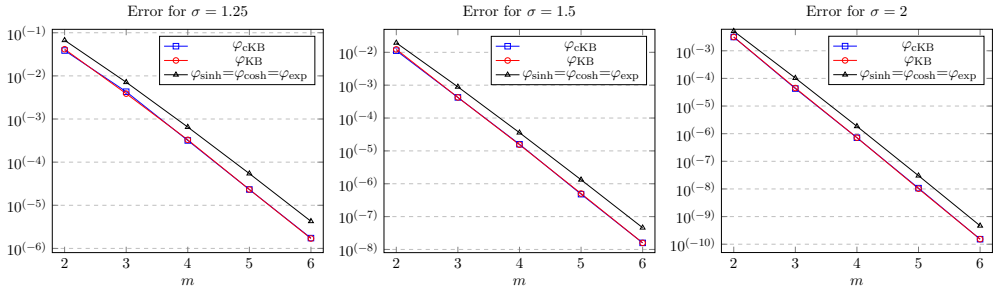
\begin{figure}[ht]
	\begin{tikzpicture}[scale=0.52]
	\begin{axis}[
	ymode=log,
	title={Error for $\sigma=1.25$},
	xlabel={$m$},
xmin=1.8, xmax=6.2,
	ymin=0.8*10^(-6), ymax=1.2*10^(-1),
	xtick={2,3,4,5,6},
	ytick={10^(-1),10^(-2),10^(-3),10^(-4), 10^(-5),10^(-6)},
	yticklabels={$10^{(-1)}$, $10^{(-2)}$, $10^{(-3)}$, $10^{(-4)}$, $10^{(-5)}$, $ 10^{(-6)}$},
	legend pos=north east,
	ymajorgrids=true,
	grid style=dashed,
	]
	\addplot[color=blue, mark=square]
	coordinates {(2, 3.8755e-02) (3,4.2966e-03) (4,3.1845e-04)  (5, 2.3238e-05)  (6,1.7344e-06)};
	\addplot[color=red, mark=o]
	coordinates {(2,4.1531e-02)  (3,3.8784e-03 )   (4,3.2948e-04)
		(5,2.3238e-05) (6,1.6964e-06)};
	\addplot[color=black, mark=triangle]
	coordinates {(2, 6.6976e-02)  (3,7.1278e-03)   (4,6.5055e-04)
		(5,5.4320e-05) (6,4.2561e-06)};
	\legend{$\varphi_{\mathrm{cKB}}$, $\varphi_{\mathrm{KB}}$, $\varphi_{\sinh}$=$\varphi_{\cosh}$=$\varphi_{\mathrm{exp}}$}
	\end{axis}
	\end{tikzpicture}
	\begin{tikzpicture}[scale=0.52]
	\begin{axis}[
	ymode=log,
	title={Error for $\sigma=1.5$},
	xlabel={$m$},
xmin=1.8, xmax=6.2,
	ymin=0.8*10^(-8), ymax=5*10^(-2),
	xtick={2,3,4,5,6},
	ytick={10^(-2),10^(-3),10^(-4),10^(-5), 10^(-6), 10^(-7), 10^(-8)},
	yticklabels={$10^{(-2)}$, $10^{(-3)}$, $10^{(-4)}$, $10^{(-5)}$, $10^{(-6)}$, $10^{(-7)}$,  $10^{(-8)}$},
	legend pos=north east,
	ymajorgrids=true,
	grid style=dashed,
	]
	\addplot[color=blue, mark=square]
	coordinates {(2,1.1239e-02) (3,4.2453e-04) (4,1.5936e-05 ) (5,4.8191e-07) (6,1.5868e-08)};
	\addplot[color=red, mark=o]
	coordinates {(2,1.2306e-02)  (3,4.2453e-04) (4,1.5444e-05 ) (5,5.0302e-07) (6,1.5868e-08)};
	\addplot[color=black, mark=triangle]
	coordinates {(2,1.9433e-02 ) (3,8.9016e-04) (4,3.5865e-05) (5,1.3237e-06) (6,4.5903e-08)};
	\legend{$\varphi_{\mathrm{cKB}}$, $\varphi_{\mathrm{KB}}$, $\varphi_{\sinh}$=$\varphi_{\cosh}$=$\varphi_{\mathrm{exp}}$}
	\end{axis}
	\end{tikzpicture}
	\begin{tikzpicture}[scale=0.52]
	\begin{axis}[
	ymode=log,
	title={Error for $\sigma=2$},
	xlabel={$m$},
xmin=1.8, xmax=6.2,
	ymin=0.5*10^(-10), ymax=6*10^(-3),
	xtick={2,3,4,5,6},
	ytick={10^(-3), 10^(-4), 10^(-5),10^(-6),10^(-7), 10^(-8), 10^(-9),10^(-10)},
	yticklabels={$10^{(-3)}$, $10^{(-4)}$,$10^{(-5)}$,$10^{(-6)}$,$10^{(-7)}$, $10^{(-8)}$,$10^{(-9)}$,$10^{(-10)}$},
	legend pos=north east,
	ymajorgrids=true,
	grid style=dashed,
	]
	\addplot[color=blue, mark=square]
	coordinates {(2, 3.1435e-03)  (3,4.2916e-05)  (4,7.1695e-07) (5,1.0699e-08 ) (6, 1.5200e-10)};
	\addplot[color=red, mark=o]
	coordinates {(2, 3.1539e-03)  (3,4.4842e-05)  (4,7.1695e-07) (5,1.0217e-08) (6, 1.5229e-10)};
	\addplot[color=black, mark=triangle]
	coordinates {(2,5.1243e-03) (3,1.0287e-04 )  (4,1.8467e-06) (5,3.0197e-08) (6,4.6553e-10)};
	\legend{$\varphi_{\mathrm{cKB}}$, $\varphi_{\mathrm{KB}}$, $\varphi_{\sinh}$=$\varphi_{\cosh}$=$\varphi_{\mathrm{exp}}$}
	\end{axis}
	\end{tikzpicture}
	
	\caption{The constants $e_{\sigma,N}({\varphi})$ of the different window functions  with shape pa\-rameter $\beta=\pi m(2-1/\sigma)$ 
for  $\sigma\in \{1.25, 1.5,2 \}$ and  $m\in \{2,\,3,\,4\}$.}
	\label{Fig.Sinhhh}
\end{figure}

In summary, the $C(\mathbb T)$-error constant of the continuous/standard Kaiser--Bessel window function is of best order ${\mathcal O}\big(m \, {\mathrm e}^{-2\pi m \sqrt{1- 1/\sigma}}\big)$. For the $\sinh$-type,
continuous/original $\exp$-type, and continuous $\cosh$-type window functions, the corresponding $C(\mathbb T)$-error constants are of order  ${\mathcal O}\big(m^{3/2} \, {\mathrm e}^{-2\pi m \sqrt{1- 1/\sigma}}\big)$. Nevertheless, our numerical results show that all window functions proposed here yield a very similar error, see Figure \ref{Fig.Sinhhh}.  Thus we can recommend the use of all these window functions.

\section*{Acknowledgments}

The authors would like to thank Stefan Kunis and Elias Wegert for several fruitful discussions on the topic.
The authors are grateful to the anonymous referees for helpful comments and suggestions.
Further, the first author acknowledges funding  by Deutsche Forschungsgemeinschaft (German Research Foundation) -- Project--ID 416228727 -- SFB 1410.

\end{document}